\documentclass{amsart}
\usepackage{amssymb}
\theoremstyle{plain}
\newtheorem{theorem}{Theorem}[section]
\newtheorem*{maintheorem}{Main Theorem}
\newtheorem{proposition}[theorem]{Proposition}
\newtheorem{lemma}[theorem]{Lemma}
\newtheorem{corollary}[theorem]{Corollary}
\theoremstyle{remark}
\newtheorem{remark}[theorem]{Remark}
\newtheorem{remarks}[theorem]{Remarks}

\theoremstyle{definition}

\newtheorem{definition}[theorem]{Definition}

\newtheorem{problem}[theorem]{Problem}
\newtheorem{problems}[theorem]{Problems}

\newtheorem{thm}[theorem]{Theorem}

\theoremstyle{definition}

\newtheorem{qst}[theorem]{Question}

\newcommand{\N}{\mathbb{N}}
\newcommand{\Q}{\mathbb{Q}}
\newcommand{\Z}{\mathbb{Z}}
\newcommand{\R}{\mathbb{R}}

\newcommand{\vp}{\varepsilon}

\newcommand{\supp}{\text{\rm supp}}

\def\ldots{\mathinner{\ldotp\ldotp\ldotp}}
\def\ldots{\mathinner{\cdotp\cdotp\cdotp}}

\def \N{\rm {\bf N}}
\def \R{\rm {\bf R}}

\def\supp{\operatorname{supp}}

\catcode`@=11 \@mparswitchfalse  
\newcounter{mnotecount}[section]

\newcommand{\mnote}[1]{} 
\title{Coefficient Quantization in Banach Spaces}
\author{S. J. Dilworth}\address{Department of Mathematics\\ University of South Carolina\\
Columbia, SC 29208\\ USA} \email{dilworth@math.sc.edu}
\author{E. Odell}
\address{Department of Mathematics \\
The University of Texas\\1 University Station C1200\\
Austin, TX 78712\\ USA}
\email{odell@math.utexas.edu}
\author{Th. Schlumprecht}
\address{Department of Mathematics, Texas A \& M University\\
College Station, TX 78712, USA}
\email{thomas.schlumprecht@math.tamu.edu}
\author{Andr\'as Zs\'ak}
\address{Fitzwilliam College, Cambridge, CB3 0DG, England}
\email{a.zsak@dpmms.cam.ac.uk}

\thanks{The research of
the second and third authors was supported by the NSF. The first, second, and fourth authors were supported by
the Linear Analysis Workshop at Texas A\&M University in 2005. All authors were supported by BIRS}
\thanks{S. J. Dilworth (Phone: 803 777 4715, Fax: 803 777 3783) is the communicating author}

\begin{document}
\begin{abstract} Let $(e_i)$ be a dictionary for a separable Banach space $X$.
We consider the problem of approximation  by  linear combinations
of dictionary elements with quantized coefficients  drawn usually from a `finite alphabet'. We investigate
several
approximation properties of this type and connect them to the Banach space geometry of $X$.  The existence of a total minimal system with
one of these properties, namely the \textit{coefficient quantization property},  is shown to be equivalent to $X$ containing $c_0$.
\end{abstract}
\maketitle

\tableofcontents
\section{Introduction}

We begin with the problem which motivates this paper.
Let $(X,\|\cdot\|)$ be a separable infinite-dimensional Banach space and let $(e_i)$
be a semi-normalized  dictionary for $X$ (i.e.\ $(e_i)$ has dense linear span in $X$).
  For a given choice of $N \in \mathbb{N}$,
consider the problem of approximating an element $x \in X$
by an element of the
`lattice'
$$\mathcal{D}^N((e_i)) = \{\sum_{i \in E} \frac{k_i}{2^N} e_i
\colon k_i \in\mathbb{Z}, \text{$E\subset\mathbb{N}$ finite}\}.$$
In many situations (e.g. when $(e_i)$ is a Schauder basis for $X$) each coefficient
$k_i/2^N$ of an approximant from $\mathcal{D}^N((e_i))$
will be  bounded by a
constant that depends only on $(e_i)$ and $\|x\|$. In this case  the approximant
will be chosen from a collection of vectors
in $\mathcal{D}^N((e_i))$
  whose coefficients
are quantized by
 a `finite alphabet'.

We investigate two natural approximation properties. The first
of these, which we call the \textit{Coefficient Quantization Property} (abbr.\ CQP),
is defined roughly as follows: for every prescribed tolerance there exists a quantization such that
every vector $x = \sum_{i\in E} a_i e_i$ in  $X$ that can be
expressed as
a finite linear combination  of dictionary elements
can be approximated
by a quantized vector $y = \sum_{i \in E} d_i e_i$ with the same (or possibly smaller) support
$E$. Thus,
for each $\varepsilon>0$, there exists $N$ such that for every $x$
with finite support $E$ there exists
$y \in \mathcal{D}^N((e_i))$  supported in $E$ such that $\|x-y\| \le \varepsilon$.

Precise definitions and some useful permanence properties
are presented in Section~\ref{sec: CQP}. One of our main results
(Theorem~\ref{thm: ballenough})  is the perhaps  surprising \mnote{Is this too strong?} fact that quantization of the unit ball
for some $\varepsilon<1$ automatically implies quantization of the whole space.

Several examples of bases  with the CQP,
including the  Schauder system for $C([0,1])$ and a class of  bases for $C(K)$, where $K$ is a countable compact metric space,
are discussed in Section~\ref{sec: examples}.
 On the other hand, it is shown that the  Haar basis for $C(\Delta)$,
where $\Delta$ denotes the Cantor set, is \textit{not} a CQP basis. It turns out that all of the natural examples
satisfy a stronger form of the CQP which we call the \textit{Strong Coefficient Quantization Property}. Roughly, this means that the quantization  of each coefficient can be an  arbitrary $\delta$-net, not necessarily a discrete subgroup of $\mathbb{R}$.

W. T. Gowers \cite{G} proved that every real-valued Lipschitz function on the
unit sphere of X is essentially constant on the sphere of an infinite-dimensional subspace of X if (and only if by \cite{OS}) X contains an isomorph
of $c_0$. A key feature of his argument was the fact that the unit vector
basis of $c_0$ has the CQP. The main results of this paper, as summarized in
the following theorem, yield an intimate connection between the CQP and
containment of $c_0$.

\begin{maintheorem} Let $X$ be a separable Banach space. Then $X$ has a fundamental and  total normalized minimal system with the CQP if and only if $c_0$ is isomorphic to a subspace of $X$. Moreover, if $X$ has a basis then $X$ has a normalized weakly null basis with the CQP if and only if $X$ contains an isomorph of $c_0$.\end{maintheorem}

The sufficiency   is proved   in Section~\ref{sec: existence}
(Theorem~\ref{thm: existence}) and the necessity is proved in Section~\ref{sec: containmentofc0}
 (Theorem~\ref{thm: conversemainthm}). The necessity result is stated more precisely  as the following dichotomy:  if $(e_i)$
is a fundamental and total minimal system with the CQP then some subsequence
of $(e_i)$ is equivalent to the unit vector basis of $c_0$ or to the summing basis of $c_0$.

For the reader who wishes to make a beeline for the proof of the Main Theorem we suggest a shorter route through the paper.
After absorbing the definitions of the CQP and SCQP in Section~2 and the NQP in Section~5, he or she should then read Section~4, Theorem~5.11 (which is very short), and Section~6.

The second natural approximation property, which we call
the \textit{Net Quantization Property} (abbr.\ NQP), is investigated in Section~\ref{sec: BQP}.
We say that $(e_i)$ has the NQP
 if for every $\varepsilon>0$
there exists $N$ such that $\mathcal{D}^N((e_i))$ is an $\varepsilon$-net for  $X$. We prove that the NQP is a weaker property than the CQP. In particular,
while the CQP is preserved under the operation of passing to a subsequence, this is not the case for the NQP. Indeed, we prove (Theorem~\ref{thm: NQPnotCQP}) that
{\it every} normalized bimonotone basic sequence may be embedded as a subsequence of a Schauder basis with the NQP.
Another  main result of Section~\ref{sec: BQP}
is related to the \textit{greedy algorithm} in Banach spaces (see e.g. \cite{DKKT}).
It is proved that the unit vector basis of $c_0$ is  the
only \textit{quasi-greedy}  NQP minimal system.

We do not know whether or not every space $X$ with an NQP basis contains $c_0$.
However, we are able to prove the weaker result that if $X$ admits a minimal system
with the NQP  then
the dual space of $X$ contains
 an isomorphic copy of $\ell_1$ (Theorem~\ref{thm: embedl1indual}). In particular, $X$ is
necessarily  non-reflexive.

 The last section contains some examples and questions of a finite-dimensional character that are related to the CQP.

Standard Banach space notation and terminology are used throughout
(see e.g.\ \cite{LT}).
For the sake of clarity, however, we recall the notation
that is used most heavily. Let $(X,\|\cdot\|)$ be a real Banach space with
 \textit{dual space} $X^*$. The \textit{unit ball} of $X$ is the
set $Ba(X)\
:=\{x\in X\colon \|x\|\le1\}.$
We write $Y \hookrightarrow X$ (where $(Y,\|\cdot\|)$ is another Banach space)
if there exists a continuous linear isomorphism from $Y$ into $X$.

Let $(e_i)$ be a
sequence in $X$. The closed linear span of   $(e_i)$
is denoted  $[(e_i)]$. We say that $(e_i)$ is \textit{weakly Cauchy}
if the scalar sequence $(x^*(e_i))$ converges for each $x^* \in X^*$. We say that $(e_i)$ is {\it nontrivial weakly Cauchy}
if $(e_i)$ is weakly Cauchy but not weakly convergent, i.e. $(e_i)$ converges weak-star to an element of $X^{**}\setminus X$.
We say that a sequence $(e_i)$
of nonzero vectors is \textit{basic}
if there exists  a positive constant $K$  such that
$$\|\sum_{i=1}^m a_i e_i \| \le K \|\sum_{i=1}^n a_ie_i\|$$
for all scalars $(a_i)$ and
all $1 \le m \le n \in \mathbb{N}$; the least such constant is called the {\it basis constant}; $(e_i)$ is \textit{monotone}
if we can take $K=1$;
$(e_i)$ is \textit{$C$-unconditional},
where $C$ is a positive constant, if
$$\|\sum_{i=1}^n \varepsilon_i a_i e_i\| \le C \|\sum_{i=1}^n a_i e_i\|$$
for all scalars $(a_i)$, all choices of signs $\varepsilon_i=\pm 1$, and all $n\ge1$.
The least such constant is called the \textit{constant of unconditionality}.
We say that $(e_i)$ is a (Schauder)  basis  for $X$ if $(e_i)$ is basic and
$[(e_i)]=X$.
Two basic sequences $(e_i)$ and $(f_i)$ are said to be \textit{equivalent}
if the mapping $e_i \mapsto f_i$ extends to a linear isomorphism from $[(e_i)]$
onto $[(f_i)]$.

For $1\le p < \infty$, $\ell_p$ is the space of real sequences
$(a_i)$ equipped with the norm  $\|(a_i)\|_p = (\sum_{i=1}^\infty |a_i|^p)^{1/p}$.
The space of sequences converging to zero (resp.\ bounded) equipped with the supremum norm
$\|\cdot\|_\infty$
is denoted $c_0$ (reps. $\ell_\infty$). The linear space of eventually zero sequences is denoted
$c_{00}$. For $(a_i) \in c_{00}$, the \textit{support}
of $x$, denoted  $\supp x$, is the set $\{i \in \mathbb{N}\colon
a_i \ne0\}$.
The space of continuous functions on a compact Hausdorff space $K$ equipped with the supremum norm $\|\cdot\|_\infty$
is denoted $C(K)$.
For Banach spaces $X$ and $Y$, the direct sum $X\oplus_\infty Y$ (resp.\ $X\oplus_1Y$) is equipped with the maximum norm
$\|(x,y)\|_\infty = \max(\|x\|,\|y\|)$ (resp.\ sum norm $\|(x,y)\|_1 = \|x\|+\|y\|$). Similarly, $(\sum_{n=1}^\infty \oplus X_n)_{0}$ and
$(\sum_{n=1}^\infty \oplus X_n)_{1}$ denote the $c_0$ and $\ell_1$ sums of the Banach spaces $(X_n)_{n=1}^\infty$ equipped with their usual norms.

Finally, it is worth emphasizing that we consider only \textit{real} Banach spaces in this paper.

\section{The Coefficient Quantization Property} \label{sec: CQP}

Throughout, $X$ will denote a separable infinite-dimensional Banach space and $(e_i)$
will denote a semi-normalized \textit{dictionary} for $X$, i.e:
\begin{itemize}
\item[(i)] there exist positive constants
$a$ and $b$ such that $a \le \|e_i\| \le b$ ($i \in \mathbb{N}$);\\
\item[(ii)] $(e_i)$ is a \textit{fundamental} system for $X$, i.e. $[(e_i)]=X$.
\end{itemize}
We say that  $(e_i)$ is a \textit{minimal system} (we shall always assume that the minimal system is semi-normalized and fundamental) if there exists a biorthogonal sequence
$(e^*_i)$ in $X^*$ such that $e^*_i(e_j) = \delta_{ij}$. We say that $(e_i)$ is \textit{total}
if $e_i^*(x) = 0$ for all $i\in\mathbb{N}$ implies that $x=0$, and that
$(e_i)$ is \textit{bounded}
if $\sup \|e_i\|\|e^*_i\| = M <\infty$. Ovsepian and Pe\l czy\'nski \cite{OP} showed that
every separable Banach space possesses a total and bounded minimal system \cite{OP}. Pe\l czy\'nski
\cite{P2} proved later that one can take $M = 1 + \varepsilon$ for any $\varepsilon>0$.

Recall that a subset $S$ of a metric space $(T,\rho)$ is a \textit{$\delta$-net}
for $A \subseteq T$ (and  is said to be \textit{$\delta$-dense} in $A$)
 if for every $x \in A$ there exists $y \in S$ such that $\rho(x,y)\le \delta.$
Also $S$ is said to be \textit{$\delta$-separated} if the distance between distinct points of $S$ is at least $\delta$.

\begin{definition} \label{def: newCQP} \mnote{This is the main definition. It is Ted's original definition except for the whole space instead of the ball}
 A dictionary $(e_i)$ has the $(\varepsilon,\delta)$-\textit{Coefficient Quantization Property}
(abbr.\ $(\varepsilon,\delta)$-CQP) if for every $x = \sum_{i \in E} a_ie_i \in X$
 (where $E$ is a finite subset of $\mathbb{N}$) there exist $n_i \in \mathbb{Z}$ ($i \in E$) such that \begin{equation} \label{eq: epsilondensenewCQP}
 \|x - \sum_{i \in E} n_i \delta e_i\| \le \varepsilon. \end{equation}
We say that $(e_i)$ has the CQP if $(e_i)$ has the $(\varepsilon, \delta)$-CQP for some $\varepsilon>0$ and $\delta>0$.
\end{definition}
\begin{remark}
Setting
$$ \mathcal F_\delta((e_i)) := \{\sum_{i \in E} n_i \delta e_i \colon \text{$E \subset \mathbb{N}$ finite, $n_i \in \mathbb{Z}$}\},$$
note that
\eqref{eq: epsilondensenewCQP} is equivalent to the following: \begin{equation*}
\text{$ \mathcal F_\delta((e_i)_{i \in E})$ is $\varepsilon$-dense in
$[(e_i)_{i\in E}]$}. \end{equation*}
\end{remark}
We begin with  some elementary observations.

\begin{proposition} \label{prop: elemproperties}
Let $(e_i)$ be  a dictionary for $X$ with the CQP and let $\varepsilon,\delta>0$.\newline
(a) The  following are equivalent:  \begin{itemize}
 \item[(i)] $(e_i)$ has  the $(\varepsilon,\delta)$-CQP. \\
\item[(ii)] $(e_i)$ has the $(\lambda \varepsilon, \lambda \delta)$-CQP for all $\lambda >0$. \\
\item[(iii)] $(e_i)$ has the $(1, \delta/\varepsilon)$-CQP. \end{itemize}
Thus, if $(e_i)$ has the CQP then there exists $c>0$ such that $(e_i)$ has the $(\varepsilon,c\varepsilon)$-CQP for
all $\varepsilon>0$. \newline
(b) The mapping
$$ \delta \mapsto \varepsilon(\delta) := \inf \{\varepsilon \colon \text{$(e_i)$ has the $(\varepsilon, \delta)$-CQP}\}.$$
is linear, i.e. $\varepsilon(\lambda \delta) = \lambda \varepsilon(\delta)$ for all $\delta>0$ and $\lambda>0$; moreover, if $(e_i)$ is linearly independent then
$(e_i)$ has the $(\varepsilon(\delta),\delta)$-CQP.
\end{proposition}
\begin{proof} (a) To prove the implication (i) $\Rightarrow$ (ii), let $\lambda>0$ and $x = \sum_{i \in E} a_i e_i$, where $E$ is finite. Since $(e_i)$ has the $(\varepsilon,\delta)$-CQP
there exist $n_i \in \mathbb{Z}$  such that $\|x/\lambda - \sum_{i \in E} n_i \delta e_i\| \le \varepsilon$.
Hence $\|x - \sum_{i \in E} n_i \lambda \delta e_i\| \le \lambda \varepsilon$, which proves (ii). The proofs of the other implications are similar.
\newline (b) The first assertion is an immediate consequence of (a), and the second  is an easy compactness argument.
\end{proof}
Now suppose that we relax Definition~\ref{def: newCQP} by only requiring that one can approximate each element $x$ of the {\it unit ball} of $X$
instead of the whole space.
Accordingly, for each $\delta>0$, we define $\varepsilon^{(b)}(\delta)$ to be the infimum of those $\varepsilon >0$ such that for all finite
 $E \subset \mathbb{N}$ we have that
$$ \text{$ \mathcal F_\delta((e_i)_{i \in E})$ is $\varepsilon$-dense in
$Ba([(e_i)_{i\in E}])$}.$$
The following theorem, which is the main result of this section, explains why  the CQP has been defined in terms of quantization of the whole
space instead of  the unit ball.
\begin{theorem} \label{thm: ballenough} \mnote{This is Thomas's theorem I have simplified the proof a little.}
Let $(e_i)$ be a dictionary for $X$. The following are equivalent: \begin{itemize}
\item[(i)] $(e_i)$ has the CQP; \\
\item[(ii)] $\varepsilon^{(b)}(\delta_0) < 1$ for some $\delta_0>0$;\\
\item[(iii)] there exists $\delta_1>0$ such that $\varepsilon(\delta) = \varepsilon^{(b)}(\delta)<\infty$ for all $0<\delta \le \delta_1$.
\end{itemize} \end{theorem}  \begin{proof}  The implications (i) $\Rightarrow$ (ii) and (iii) $\Rightarrow$ (i) are clear.
To prove the nontrivial implication  (ii) $\Rightarrow$ (iii), let
 $q_0:= (\vp^{(b)}(\delta_0)+1)/2<1$.
First we show that there  exist $0<q_1<1$ and  $\delta_1>0$ such that   for
every $0<\delta<\delta_1$, we have
$\vp^{(b)}(\delta) < q_1$.

Indeed,  choose $n_1\in\mathbb{N}$ and $0<q_1<1$ such that $$\frac{n_1+1}{n_1} q_0<q_1<1,$$ and set  $\delta_1:=\dfrac{\delta_0}{n_1}$.
For $0<\delta\le \delta_1$ and $x=\sum_{i\in E}a_ie_i\in Ba(X)$, with $E\subset \mathbb{N}$  finite,
choose  $n\in\mathbb{N}$ such that $\dfrac{\delta_0}{n+1}<\delta\le \dfrac{\delta_0}{n}$ (note that $n \ge n_1$) and choose
 $k_i\in \mathbb{Z}$ ($i \in E$) such that
$$\Bigl\|\sum_{i\in E} a_i e_i\frac{\delta_0}{(n+1)\delta} -
\sum_{i\in E} k_i\delta_0 e_i\Bigr\|<q_0.$$
Thus, since $n\ge n_1$,
$$\Bigl\|\sum_{i\in E} a_i e_i-\sum_{i\in E} k_i(n+1)\delta e_i\Bigr\|\le q_0\frac{(n+1)\delta}{\delta_0}
\le q_0\frac{n+1}n < q_1,$$
which implies that $\varepsilon^{(b)}(\delta) < q_1$.

Suppose that $0<\delta,\tilde\delta \le\delta_1$ satisfy \begin{equation} \label{eq: deltacondition}
q_1 \le \frac{\delta}{\tilde\delta }\le \frac{1}{q_1}. \end{equation}
We claim that \begin{equation} \label{eq: deltaeq}
\frac{\vp^{(b)}(\delta)}{\delta}\le \frac{\vp^{(b)}(\tilde\delta )}{\tilde\delta }. \end{equation}
Once the claim is shown, it follows, by exchanging  the roles of
$\delta$ and $\tilde\delta $, that we  also have
$$\frac{\vp^{(b)}(\tilde\delta )}{\tilde\delta }\le \frac{\vp^{(b)}(\delta)}{\delta},$$
which implies  local linearity and, thus, linearity of $\vp^{(b)}$ on $(0,\delta_1]$.

Let $x = \sum_{i \in E} a_ie_i \in Ba(X)$ with $E$ finite. There exists $y = \sum_{i \in E} k_i \delta e_i\in \mathcal{F}_{\delta}((e_i))$
such that $\|x - y\| < q_1$.  Note that $(\tilde \delta/\delta)(x-y) \in Ba(X)$ by \eqref{eq: deltacondition}.
Hence,
given $\eta>0$,  there exists $z=\sum_{i \in E} m_i \tilde\delta e_i \in \mathcal{F}_{\tilde\delta}((e_i))$ such that
$$\| \frac{\tilde \delta}{\delta}(x-y) - z\| < (1+\eta)\varepsilon^{(b)}(\tilde \delta),$$
i.e.
$$\| x - \sum_{i \in E} (k_i+m_i) \delta e_i\| < (1+ \eta) \frac{\delta}{\tilde \delta} \varepsilon^{(b)}(\tilde \delta),$$
which yields \eqref{eq: deltaeq} since $\eta>0$ is arbitrary.

In order show that $\vp(\cdot)=\vp^{(b)}(\cdot)$ on $(0,\delta_1]$, let
  $0<\delta\le \delta_1$,  let
$x=\sum_{i\in E} a_i e_i$, with $E\subset \mathbb{N}$  finite, and let $\eta>0$ be arbitrary.
  If $\| x\|\ge 1$ there exist $k_i\in\mathbb{Z}$ ($i\in E$) such that
$$\Big\|\frac{x}{\|x\|}-\sum_{i\in E} k_i\frac{\delta}{\|x\|}e_i\Big\|< (1+\eta)\vp^{(b)}\Big(
\frac{\delta}{\|x\|}\Big)
=(1+\eta)\frac{\vp^{(b)}(\delta)}{\|x\|}$$
and thus
\begin{equation}\label{E:1}\Big\|x-\sum_{i\in E} k_i{\delta}e_i\Big\|\le (1+\eta)\vp^{(b)}(\delta).
\end{equation}
If $\|x\|\le 1$ we can of course also find  $k_i \in \mathbb{Z}$ such that
(\ref{E:1}) holds. Since $\eta>0$ is arbitrary,
it follows that  $\vp(\cdot)\le \vp^{(b)}(\cdot)$ and, thus, $\vp(\cdot)= \vp^{(b)}(\cdot)$
on $(0,\delta_1]$. \end{proof}
The following corollary is a quantitative version of the last result.
\begin{corollary} \label{cor: preciseversion} \mnote{Here is a quantitative version}
Let $0<\varepsilon_0<1$ and $\delta>0$. If $\mathcal{F}_\delta((e_i)_{i \in E})$ is $\varepsilon_0$-dense in $Ba([(e_i)_{i \in E}])$
for all finite $E \subset \mathbb{N}$
then   $\mathcal{F}_\delta((e_i)_{i \in E})$ is $\varepsilon_1$-dense in $[(e_i)_{i\in E}]$ for all \begin{equation} \label{eq: epsilon1}
\varepsilon_1> (\left\lfloor \frac{\varepsilon_0}{1-\varepsilon_0}\right\rfloor+1)\varepsilon_0. \end{equation}
(Here $\lfloor x\rfloor$ denotes the integer part of $x$.)
In particular, if $\varepsilon_0<1/2$, then $\mathcal{F}_\delta((e_i)_{i \in E})$ is $\varepsilon_1$-dense in $[(e_i)_{i \in E}]$ for all $\varepsilon_1>\varepsilon_0$.
\end{corollary}
\begin{proof} Using the notation of the last proof, we may take \linebreak $n_1 = \lfloor\varepsilon_0/(1-\varepsilon_0)\rfloor +1$. The last proof  yields $$\varepsilon(\delta/n_1) =
\varepsilon^{(b)}(\delta/n_1) \le \varepsilon^{(b)}(\delta) \le \varepsilon_0.$$ Thus, $\varepsilon(\delta) \le n_1 \varepsilon_0$, which gives the result.
\end{proof}
\begin{remark}  \label{rem: generallattice}  The assumption that $(e_i)$ is semi-normalized is not required
for the validity of Corollary~\ref{cor: preciseversion}. Moreover, if $(e_i)$ is linearly independent then  strict inequality in
\eqref{eq: epsilon1} may be replaced by non-strict inequality. Finally, the result is also valid for quasi-normed spaces.
\end{remark}
In the finite-dimensional setting Corollary \ref{cor: preciseversion}   can be formulated as a covering result of independent interest.
\begin{theorem} \mnote{I thought this was interesting.}
Let $K \subset \mathbb{R}^n$ be a compact zero-neighborhood that is star-shaped about zero (i.e. $\lambda K \subseteq K$ for all $0 \le \lambda \le 1$)
and let $L \subset \mathbb{R}^n$ be a lattice (i.e.\ a discrete subgroup of $\mathbb{R}^n$). If $K \subset L+ \varepsilon_0 K$, where $0 < \varepsilon_0<1$, then $\mathbb{R}^n = L + \varepsilon_1 K$,
where $\varepsilon_1=(\left\lfloor \varepsilon_0/(1-\varepsilon_0)\right\rfloor+1)\varepsilon_0$.
\end{theorem}
\begin{proof} The gauge functional $\|x\|_K := \min \{t>0 \colon x \in tK\}$ is positively homogeneous, which is the only property of the norm that is used in the proof of Theorem~\ref{thm: ballenough}. Hence, setting
$$\varepsilon^{(b)}_L(\delta) := \min\{\varepsilon \colon K \subset \delta L + \varepsilon K\},$$the proof of Theorem~\ref{thm: ballenough} yields $$\varepsilon^{(b)}_L(\delta) = n_1\delta \varepsilon^{(b)}_L(1/n_1)\le n_1 \delta \varepsilon_0$$ for all $0 \le \delta \le 1/n_1$,
where $n_1 := n_1(\varepsilon_0)$ is defined as in the proof of Corollary~\ref{cor: preciseversion}.
The proof is concluded as before.
\end{proof}
The examples presented in the next section all have a formally stronger version of the CQP which we now define.
\begin{definition} \mnote{This is the "strong" version, like my original definition except now for the whole space.}
 Let $\varepsilon>0$ and let $\delta>0$. \newline
(a) A dictionary $(e_i)$ has the $(\varepsilon,\delta)$-\textit{Strong Coefficient Quantization Property}
(abbr.\ $(\varepsilon,\delta)$-SCQP) if for every sequence $\overline{D} :=(D_i)$ of $\delta$-nets for $\mathbb{R}$, such that $0 \in D_i$, and for every $x = \sum_{i \in E} a_ie_i$ in
$X$ (where $E$ is a finite subset of $\mathbb{N}$) there exist $d_i \in D_i$ ($i \in E$) such that \begin{equation} \label{eq: epsilondense}
 \|x - \sum_{i \in E} d_i e_i\| \le \varepsilon. \end{equation}
(b)  $(e_i)$ has the SCQP if $(e_i)$ has the $(\varepsilon,\delta)$-SCQP for some $\varepsilon>0$ and $\delta>0$.
\end{definition}

\begin{remarks} \label{rems: elemremarks}
(i) If we set
$$ \mathcal F_{\overline{D}}((e_i)) := \{\sum_{i \in E} d_i e_i \colon \text{$E\subset\mathbb{N}$ finite}, d_i \in D_i\},$$
then
\eqref{eq: epsilondense} is equivalent to the following:\begin{equation*}
\text{$ \mathcal F_{\overline{D}}((e_i)_{i \in E})$ is $\varepsilon$-dense in
$[(e_i)_{i\in E}]$}. \end{equation*} \newline
(ii) The obvious analogue for the SCQP of   Proposition~\ref{prop: elemproperties} is valid. \newline
 (iii)
 Note also the implication
$(\varepsilon, \delta)$-SCQP $\Rightarrow (\varepsilon, 2\delta)$-CQP since $2\delta \mathbb{Z}$ is a $\delta$-net. \newline
(iv) If $(e_i)$ has the $(\varepsilon, \delta)$-CQP, we  say that
$(e_i)$ is an \textit{$(\varepsilon,\delta)$-CQP dictionary}, and similarly for the SCQP. \newline
(v) To avoid  repetition we shall assume henceforth \textit{that every $\delta$-net for $\mathbb{R}$ contains zero}. \newline
(vi) Unless stated otherwise all sums of the form $\sum a_i e_i$ will be assumed to be
\textit{finite}.
\end{remarks} The uniformity built  into the definition of the SCQP (i.e.  that $\varepsilon$ depends only on $\delta$, not on the choice of  $(D_i)$) is natural in view of the following uniform boundedness  result.
\begin{proposition} \label{prop: uniformity} \mnote{This is adapted from Thomas's notes}
 Let $(e_i)$ be a dictionary for $X$. The following are  equivalent:
\begin{itemize}
\item[(i)] $(e_i)$ has the SCQP;
\item[(ii)] For all $\delta>0$ and for every sequence  $(D_i)$ of $\delta$-nets there exists $M>0$ such that  for every $x = \sum_{i \in E} a_ie_i \in
X$ (where $E$ is a finite subset of $\mathbb{N}$) there exist $d_i \in D_i$ ($i \in E$) such that \begin{equation*}
 \|x - \sum_{i \in E} d_i e_i\| \le M;\end{equation*}
\item[(iii)] Condition (ii) for $\delta=1$.
\end{itemize} \end{proposition}
\begin{proof} Clearly, (i) $\Rightarrow$ (ii) $\Rightarrow$ (iii). To prove (iii) $\Rightarrow$ (i), we argue by contradiction. Suppose that (i)
does not hold. Then by (ii) of Remarks~\ref{rems: elemremarks}
 $(e_i)$ fails the $(M,1)$-SCQP for all $M>0$. First we  construct by induction a sequence $(E_n)$ of  finite disjoint subsets of $\mathbb{N}$, a sequence $((D^n_i))$ of sequences of $1$-nets,
and vectors $x_n = \sum_{i \in E_n} a^n_i e_i \in X$ ($n \ge1$) such that \begin{equation} \label{eq: badsequence}
\inf \{ \|x_n - \sum_{i \in E_n} d^n_i e_i\| \colon d^n_i \in D^n_i \}>n \qquad (n\ge1). \end{equation}
Suppose that $n_0 \ge1$ and that the construction has been carried out for all $n < n_0$. Set $F := \cup_{n<n_0} E_n$. Since $(e_i)$ does not have the
$(M,1)$-SCQP for $M = \operatorname{card}(F) \max_{i \in F} \|e_i\| + n_0$
there exist a sequence $(D^{n_0}_i)$ of $1$-nets, a finite set $G \subset \mathbb{N}$, and $x = \sum_{i \in G} a_i e_i \in X$ such that
\begin{equation} \label{eq: Gapprox}
\inf \{ \|x - \sum_{i \in G} d^{n_0}_i e_i\| \colon d^{n_0}_i \in D^{n_0}_i \}> \operatorname{card}(F) \max_{i \in F} \|e_i\| + n_0. \end{equation}
Choose $d^{n_0}_i \in D^{n_0}_i$ such that $|a_i - d^{n_0}_i| \le1$ for $i \in G \cap F$.  Then
$$\|\sum_{i \in G \cap F} (a_i -d^{n_0}_i) e_i\| \le \operatorname{card} (G \cap F) \max_{i \in G \cap F} \|e_i\|,$$
and thus \eqref{eq: Gapprox} yields
$$\inf \{ \|\sum_{i \in G \setminus F}(a_i- d^{n_0}_i) e_i\| \colon d^{n_0}_i \in D^{n_0}_i \}>  n_0.$$
Set $E_{n_0} := G \setminus F$ and $x_{n_0} = \sum_{i \in G \setminus F} a_i e_i$ to complete the induction.
Now define a sequence $(D_i)$ of $1$-nets  as follows: \begin{equation*}
D_i = \begin{cases} D^n_i &\text{if there exist $n$ such that $i \in E_n$},\\
2\mathbb{Z} &\text{otherwise.} \end{cases} \end{equation*}
Then by \eqref{eq: badsequence} $(D_i)$ does not satisfy (iii).

\end{proof}

Our first permanence result ensures that the SCQP is preserved under linear isomorphisms.

\begin{proposition} \label{prop: isom} Suppose that $T : X \rightarrow Y$ is a bounded operator. Suppose also that $(e_i)$ is a dictionary for $X$ with
the property   that $(T(e_i))$ is a  dictionary for $Y$.
\newline
(a) If
$(e_i)$ is an $(\varepsilon, \delta)$-SCQP dictionary for $X$ then
$(T(e_i))$ is an $(\varepsilon\|T\|,\delta)$-SCQP dictionary for $Y$. \newline
(b) If $(e_i)$ has the SCQP then $(T(e_i))$ also has the SCQP. \end{proposition}
\begin{proof} (a) Let $(D_i)$ be any family of $\delta$-nets for $\mathbb{R}$.
Consider  $\sum_{i \in E} a_iT(e_i) \in Y$, where $E$ is a finite subset of $\mathbb{N}$.
 Since $(e_i)$ has the
$(\varepsilon, \delta)$-CQP there exist $d_i \in D_i$ such that
$$\|\sum_{i \in E} (a_i -d_i)e_i \| \le \varepsilon,$$
whence
$$ \|\sum_{i \in E} (a_i-d_i)T(e_i)\|\le \|T\|\varepsilon.$$  \newline
(b) This follows at once from (a).
\end{proof}

\begin{remark} The analogue of Proposition~\ref{prop: isom} for the CQP is also valid. \end{remark}

The following useful result shows that the SCQP is preserved after normalization of the dictionary.

\begin{proposition} \label{prop: normalize}
Suppose that $(e_i)$ has the $(\varepsilon,\delta)$-SCQP and that $a \le \|e_i\| \le b$.
Then the normalized dictionary $(e_i/\|e_i\|)$ has the $(\varepsilon,\delta')$-SCQP for $\delta'=a\delta$.
\end{proposition}

\begin{proof} Let $(D_i')$ be a family of $\delta'$-nets for $\mathbb{R}$. Then each
$D_i=\{d_i'/\|e_i\| \colon d'_i \in D'_i\}$ is a $\delta$-net. Since $(e_i)$ has the
$(\varepsilon,\delta)$-SCQP, it follows that for each
$\sum_{i \in E} a_i (e_i/\|e_i\|)$ in $X$, where $E$ is a finite subset of $\mathbb{N}$, there exist $d'_i \in D'_i$ ($i \in E$)
such that
$$ \| \sum_{i \in E}\frac{a_i}{\|e_i\|}e_i - \sum_{i \in E} \frac{d'_i}{\|e_i\|} e_i\|
\le \varepsilon.$$ \end{proof}
We conclude this section with some open problems. \mnote{I gathered together some problems here that I don't think we know the answer to?}
\begin{problems} (1) For a given dictionary $(e_i)$ is the SCQP equivalent to the CQP? \newline
(2) Does the analogue of Theorem~\ref{thm: ballenough} for the SCQP hold? \newline
(3) Does the analogue
of Proposition~\ref{prop: normalize} for the CQP  hold?
\end{problems}
\begin{remark} We say that a dictionary $(e_i)$ has property P if the following condition holds. There exists $\delta>0$ such that for all
$\delta$-nets $(D_i)$ and for all finite  $E \subseteq \mathbb{N}$ there exist $d_i \in D_i \setminus \{0\}$ ($i \in E$) such that
$\|\sum_{i \in E} d_i e_i\|\le1$. To see that Property P implies the SCQP, let $(D_i)$ be a sequence of $\delta$-nets and consider $x = \sum_{i \in E} a_i e_i$.
Clearly, each  $D'_i := \{d_i - a_i\colon d_i \in D_i\} \cup \{0\}$ is a $\delta$-net.  Property P implies that there exist $d_i \in D_i$ ($i \in E$) with $d_i \ne a_i$
such that $\|\sum_{i \in E} (d_i - a_i) e_i \| \le 1$, so $(e_i)$ has the SCQP. When $(e_i)$ is linearly independent, one can also show that the converse implication holds, i.e. that the SCQP implies Property P. So for a linearly independent dictionary the first problem stated above is equivalent to the following: is the CQP equivalent to Property P?
\end{remark}
\section{Examples} \label{sec: examples}
\subsection{The unit vector basis of $c_0$}
The unit vector basis of $c_0$ has the $(\varepsilon,\varepsilon)$-SCQP.
To see this, let $(D_i)$ be a sequence of $\varepsilon$-nets. Given
 $x = \sum_{i \in E} a_i e_i$, simply
choose $d_i \in D_i$   such that $|a_i - d_i| \le \varepsilon$. Then \begin{equation*}
\|\sum_{i \in E} a_i e_i - \sum_{i \in E} d_ie_i\| = \max_{i \in E} |a_i-d_i| \le\varepsilon.
\end{equation*}
It is instructive to note that if $(e_i)$ is a bounded minimal system then  the above procedure for choosing the approximation
 is only effective  for the unit vector basis of
$c_0$. To be precise, suppose that the $\delta$-nets $(D_i)$ are
$\gamma$-separated for some $\gamma>0$. Consider the following algorithm:   choose $d_i$
 to be the  \textit{best approximation}
 to the coefficient $a_i$   (or the best approximation of smallest absolute value when $a_i$ is exactly half-way between two
$d_i$ values).

\begin{proposition} Let $(e_i)$ be a bounded  minimal system. The following are equivalent:
\begin{itemize}
\item[(i)] $(e_i)$ is equivalent to the unit vector basis of $c_0$;
\item[(ii)] $(e_i)$ has the SCQP and the algorithm described above implements the SCQP (when the $\delta$-nets are $\gamma$-separated);
\item[(iii)] $(e_i)$ has the CQP and the algorithm described above implements the CQP (for  $D_i = \mathbb{Z}\delta$).
\end{itemize} \end{proposition}
\begin{proof} $(i)\Rightarrow (ii)$ was proved above and $(ii) \Rightarrow (iii)$ is trivial. For the   proof of
$(iii) \Rightarrow (i)$, suppose that the $(\varepsilon,\delta)$-CQP for $(e_i)$ is implemented by the aforementioned algorithm, where
$0<\varepsilon<1$ and $\delta>0$.
Let $ x=\sum_{i \in E} a_i e_i$ be a unit vector and suppose that $\max |a_i| < \delta/2$.
According to the algorithm, we should
approximate $x$ by taking $d_i =0$ for all $i \in E$, which  yields the contradiction
$1=\|x\| \le \varepsilon<1$. Hence \begin{equation*}
\frac{1}{M} \max |a_i| \le \|x\| \le \frac{2}{\delta} \max |a_i|, \end{equation*}
where $M = \sup \|e_i^*\|$.
Thus, $(i)$ holds. \end{proof}
\subsection{The summing basis of $c_0$} The linear space of sequences $(a_i)$ for which  $\sum_{i=1}^\infty a_i$
converges is a Banach space when equipped with the following norm:
$$\|(a_i)\|_{sb} = \sup_n |\sum_{i=1}^n a_i|.$$  This space is isometrically isomorphic to the space $c$ of convergent sequences with the supremum norm. The unit vector basis $(e_i)$
 is equivalent to a conditional basis of $c_0$ called the summing basis.

To see that $(e_i)$ has the $(\varepsilon,\varepsilon)$-SCQP, let $(a_i) \in c_{00}$. Suppose that $(d_i)_{i=1}^k$ have been chosen so that
$|\sum_{i=1}^j(a_i-d_i)| \le \varepsilon$ for $1 \le j \le k$. Then we continue by choosing $d_{k+1}\in D_{k+1}$ so that $d_{k+1}=0$ if $a_{k+1}=0$ and so that
$|\sum_{i=1}^{k+1}(a_i - d_i)| \le \varepsilon$.

Let us generalize this example as follows. Let $N \in \mathbb{N}$. For each $1 \le n \le N$, let $(\varepsilon^n_i)_{i=1}^\infty$ be a sequence
of signs $\varepsilon^n_i = \pm1$. Consider the following norm on $c_{00}$:
$$\|(a_i)\| = \max_{1 \le n \le N} \|(\varepsilon^n_i a_i)\|_{sb}.$$
For each $\eta = (\eta_n)_{n=1}^N \in \{-1,1\}^N$, let $A_\eta = \{m \in \mathbb{N} \colon \varepsilon^n_m = \eta_n, 1 \le n \le N\}$.
Then $(A_\eta)$ ($\eta \in \{1,-1\}^N$) is a partition of $\mathbb{N}$. Note that for $(a_i), (d_i) \in c_{00}$, the triangle inequality gives
\begin{equation} \label{eq: generalized}
\|(a_i-d_i)\| \le \sum_{\eta \in \{1,-1\}^N} \|(a_i - d_i)_{i \in A_\eta}\|_{sb}
\end{equation}
Now suppose that $(D_i)$ is a sequence of $\varepsilon/2^N$-nets for $\mathbb{R}$.
For each $\eta \in \{1,-1\}^N$, choose $d_i \in D_i$ for $i \in A_\eta$ so that $\|(a_i - d_i)_{i \in A_\eta}\|_{sb}\le \varepsilon/2^N$.
This is possible since the summing basis has the $(\varepsilon/2^N,\varepsilon/2^N)$-SCQP.
It follows from \eqref{eq: generalized}  that $\|(a_i-d_i)\| \le \varepsilon$. Hence $\|\cdot\|$ has the $(\varepsilon,\varepsilon/2^N)$-SCQP.
\subsection{The Schauder basis}

Let us recall the definition of the classical Schauder basis $(f_i)_{i\ge0}$ for $C([0,1])$:
$f_0(t)=1$, $f_1(t)=t$, and for $i=2^k + l, 0 \le l < 2^k$,  $f_i$ is the piecewise-linear function
supported  on $[l2^{-k},(l+1)2^{-k}]$ satisfying
$f_i(l2^{-k})=f_i((l+1)2^{-k})=0$ and $f_i((2l+1)2^{-k-1})=1$.

\begin{theorem} The Schauder basis for $C([0,1])$ has the $(\varepsilon,\varepsilon)$-SCQP for all $\varepsilon>0$.
\end{theorem}
\begin{proof} Let $(D_i)$ be a sequence of $\varepsilon$-nets. Suppose that $N \ge 0$ and
that  $x = \sum_{i=0}^N a_if_i$. We shall prove that there
exist $d_i \in D_i$ such that \begin{equation} \label{eqn: Schauder1}
\|\sum_{i=0}^k(a_i-d_i) f_i\|_\infty \le \varepsilon
\end{equation}
for $0 \le k \le N$ and such that $d_i=0$ if $a_i=0$. Choose $d_0 \in D_0$
such that $|a_0-d_0| \le \varepsilon$ and choose $d_1 \in D_1$ such that
$|a_0+a_1-d_0-d_1| \le \varepsilon$ (with $d_i =0$ if $a_i=0$). This establishes \eqref{eqn: Schauder1}
for $k=0$ and $k=1$. Suppose that $2 \le n \le N$ and that $d_0,\dots,d_{n-1}$ have been chosen so that
\eqref{eqn: Schauder1} holds for $0\le k \le n-1$. Let the support of $f_n$ be the dyadic interval $[a,b]$ and consider
\begin{equation*}
g(x)=|\sum_{i=0}^n a_i f_i(x) - \sum_{i=0}^{n-1}d_i f_i(x)|. \end{equation*}
Then $g$ is piecewise-linear on $[a,b]$ with nodes at $a$, $b$, and $(a+b)/2$. So $g$ must attain its maximum at one of these three points. If the maximum occurs at either $x=a$ or $x=b$, then,
since $f_n(a)=f_n(b)=0$, it follows from the case $k=n-1$ of \eqref{eqn: Schauder1}
that  \begin{equation*}
\max_{x \in [a,b]} g(x) \le
\max_{x \in [0,1]}\|\sum_{i=0}^{n-1}(a_i-d_i) f_i\|_\infty \le \varepsilon.
\end{equation*} Then, setting $d_n=0$, \eqref{eqn: Schauder1} will be satisfied for
$k=n$.
So suppose that the maximum is attained at $(a+b)/2$. Choose $d_n \in D_n$
such that \begin{equation*}
|\sum_{i=0}^{n-1} (a_i-d_i) f_i(\frac{a+b}{2}) +a_n f_n(\frac{a+b}{2}) - d_n| \le \varepsilon.
\end{equation*} With this choice of $d_n$, we see that \eqref{eqn: Schauder1}
is again satisfied for $k=n$.
\end{proof} \begin{remark} Let $K$ be an uncountable compact metric space. Then $C(K)$ is uniformly isomorphic to $C([0,1])$ by
Milutin's Theorem \cite{M}. Since the Schauder basis of $C([0,1])$ has the $(\varepsilon, \varepsilon)$-SCQP,
it follows from Propositions~\ref{prop: isom} and ~\ref{prop: normalize}
that $C(K)$ has a normalized $(\varepsilon, c\varepsilon)$-SCQP basis for some absolute constant $c>0$.
\end{remark}
\subsection{Tree spaces}
By a \textit{tree} we shall mean a partially ordered set $(\mathcal{T},\le)$
with the property that each \textit{node} $\alpha \in \mathcal{T}$ has  finitely many linearly ordered predecessors (with respect to $\le$).
We say that $\mathcal{T}$ is \textit{rooted} if there
is exactly one  node without an immediate predecessor.
The  tree $\mathcal{T}_\infty$ is the rooted tree with the property that every node has countably infinitely many
immediate successors.
  We equip $c_{00}(\mathcal{T})$ with the following
norm: \begin{equation*}
\|x\| = \max_{\beta \in \mathcal{T}}|\mathcal{S}_\beta(x)|, \end{equation*}
where $\mathcal{S}_\beta(x)=\sum_{\alpha\le\beta}x(\alpha)$. Let $\mathcal{S}(\mathcal{T})$ denote
the completion of the normed space $(c_{00}(\mathcal{T}),\|\cdot\|)$.

Henceforth we shall assume that $\mathcal{T}$ is countably infinite. Suppose that
 $(\alpha(i))$ is any enumeration of $\mathcal{T}$
which respects the ordering of $\mathcal{T}$, i.e. such that $$\alpha(i) \le \alpha(j) \Rightarrow i \le j.$$
Clearly, $(e_{\alpha(i)})$
is  a normalized monotone basis for $\mathcal{S}(\mathcal{T})$.
\begin{proposition} (a) Suppose that $\mathcal{T}$ is rooted. Then
 $\mathcal{S}(\mathcal{T})$ is isometrically isomorphic to $C(K)$, where $K$ is the weak-star closure of $\{\mathcal{S}_\beta \colon \beta \in \mathcal{T}\}$ in $Ba(\mathcal{S}(\mathcal{T})^*)$.  \newline
(b) If $K$ is a countable compact metric space then  $C(K)$ is isometrically isomorphic to $\mathcal{S}(\mathcal{T})$ for some rooted tree $\mathcal{T}$. \newline
(c)  $\mathcal{S}(\mathcal{T}_\infty)$ is isometrically isomorphic to
$C(\Delta)$, where $\Delta$ denotes the Cantor set.\end{proposition} \begin{proof} (a) It is easily seen that $c_{00}(\mathcal{T})$ is a separating subalgebra of $C(K)$.
Since $\mathcal{S}_\alpha(e_\emptyset) = 1$ for all $\alpha \in \mathcal{T}$,
where $\emptyset$ is the root node, it follows that $\chi_K \in  c_{00}(\mathcal{T})$, and
hence  by the Stone-Weierstra$\ss$ theorem that $c_{00}(\mathcal{T})$ is dense in $C(K)$.
\newline
(b) \mnote{I have added this proof.}
It is well-known that every countable compact metric space is homeomorphic to an ordinal
interval $[0,\alpha]$,  for some countable ordinal $\alpha$, with the order topology. We prove the result by transfinite induction.
Suppose the result holds for  $K=[0,\beta]$ for all $0 \le \beta<\alpha$. There exist $1 \le n \le \infty$ and countable ordinals $\alpha_j<\alpha$ ($0 \le j <n$) such that $K:=[0,\alpha]$ is homeomorphic to the one-point compactificiation of the disjoint union of
the ordinal intervals $K_j := [0, \alpha_j]$ ($0 \le j < n$). By hypothesis there exist trees $\mathcal{T}_j$ ($0 \le j <n$) such that $\mathcal{S}(\mathcal{T}_j)$
is isometrically isomorphic to $C(K_j)$. Let $\mathcal{T}$ be the rooted tree which has each $\mathcal{T}_j$ ($0 \le j <n$) as a subtree immediately succeeding the root node.
Then $\mathcal{S}(\mathcal{T})$ is easily seen to be isometrically isomorphic to $C(K)$.
 \newline
(c) In this case $K$ is easily seen to be a perfect and totally disconnected compact metric space,  and thus homeomorphic to $\Delta$. \end{proof}
\begin{theorem} $(e_\alpha)_{\alpha \in \mathcal{T}}$ has the $(\varepsilon,\varepsilon)-SCQP$ in $\mathcal{S}(\mathcal{T})$ for all $\varepsilon>0$. \end{theorem}
\begin{proof} Let $(\alpha(i))$ be any ordering of the basis which respects the ordering
of $\mathcal{T}$.
Let $\varepsilon>0$ and let $(D_\alpha)_{\alpha\in \mathcal{T}}$ be a family of $\varepsilon$-nets and
suppose that
$\sum_{i \in E} x_{\alpha(i)} \in c_{00}(\mathcal{T})$. We define $d_\alpha \in D_\alpha$ inductively.
Suppose that $n\ge0$ and that $d_{\alpha(1)},\dots,d_{\alpha(n)}$ have been chosen such that
\begin{equation} \label{eq: treebasis}
|\mathcal{S}_\gamma(\sum_{i=1}^n(x_{\alpha(i)}-d_{\alpha(i)}) e_{\alpha(i)})| \le \varepsilon
\end{equation}
for all $\gamma \in \mathcal{T}$ (This condition is vacuous for $n=0$.) If $x_{\alpha(n+1)}=0$, set
$d_{\alpha(n+1)}=0$. Otherwise choose $d_{\alpha(n+1)} \in D_{\alpha(n+1)}$ such that
$$ |\sum_{\beta<\alpha(n+1)}(x_\beta-d_\beta) + x_{\alpha(n+1)}-d_{\alpha(n+1)}|\le\varepsilon,$$
noting that if $\beta<\alpha(n+1)$ then $\beta=\alpha(j)$ for some $j \le n$. Now we verify the inductive
hypothesis for $n+1$. If $\gamma \ge \alpha(n+1)$ then
$$ |\mathcal{S}_\gamma(\sum_{i=1}^{n+1}(x_{\alpha(i)}-d_{\alpha(i)}e_{\alpha(i)}))|
=|\sum_{\beta\le\alpha(n+1)}(x_\beta-d_\beta)| \le\varepsilon.$$
On the other hand, if $\gamma < \alpha(n+1)$ or if $\gamma$ and $\alpha(n+1)$ are incomparable, then
$$|\mathcal{S}_\gamma(\sum_{i=1}^{n+1}(x_{\alpha(i)}-d_{\alpha(i)}) e_{\alpha(i)})|=
|\mathcal{S}_\gamma(\sum_{i=1}^n(x_{\alpha(i)}-d_{\alpha(i)}) e_{\alpha(i)})|\le\varepsilon$$
by the inductive assumption \eqref{eq: treebasis}. This completes the verification of the inductive step. It follows that
$$
\|\sum_{n \in E} (x_{\alpha(n)}-d_{\alpha(n)})e_{\alpha(n)}\|\le\varepsilon.$$ Thus,
$(e_\alpha)_{\alpha \in \mathcal{T}}$ has the $(\varepsilon,\varepsilon)$-SCQP.
\end{proof}
\begin{corollary} If $K$ is a countable compact metric space or if $K=\Delta$ then $C(K)$ has a monotone basis with the
$(\varepsilon,\varepsilon)$-CQP for all $\varepsilon>0$. \end{corollary}
\begin{remark} \label{rem: neighborly} In all of the above examples the dictionary $(e_i)$ has the {\it neighborly CQP}, i.e. for every
$x = \sum_{i \in E} a_ie_i$  with finite support, the approximation $y = \sum_{i \in E} n_i\delta e_i$ satisfies $|a_i - n_i\delta|\le \delta$.
We do not know whether  this holds in general, i.e. whether the CQP implies  the neighborly CQP. \end{remark}
\subsection{The Haar Basis for $C(\Delta)$} We have already seen that $C(\Delta)$  has a monotone basis with the $(\varepsilon,\varepsilon)$-CQP. Surprisingly, however, the natural basis of $C(\Delta))$, namely the Haar basis,  does \textit{not} have the CQP. Let us recall the definition of the Haar basis.
Let $\Delta_0:=\Delta$, and,  for $k\ge0$, let $\Delta_{2k+1}$ and $\Delta_{2k+2}$ be the left-hand and
right-hand halves of $\Delta_k$ obtained by removing the `middle third' in the
classical construction of the Cantor set.  Then
\begin{equation*}
h_i= \begin{cases} \chi_{\Delta}&\text{for i=0}\\
\chi_{\Delta_{2i-1}}-\chi_{\Delta_{2i}} &\text{for $i>0$}. \end{cases}
\end{equation*}
Clearly, $(h_i)_{i=0}^\infty$ is
a monotone basis
for $C(\Delta)$.  For $k=1,2,\dots$, we say that the $2^{k-1}$ Haar functions $\{h_i \colon 2^{k-1} \le i < 2^k\}$
are on the $k$-th level.

\begin{proposition} \label{prop: haarsystem} Let $0<\varepsilon<1$ and let $\delta>0$. Then $\mathcal{F}_\delta((h_i))$
is not an $\varepsilon$-net for the unit ball of $C(\Delta)$. In particular, $(h_i)$ does not have the CQP.\end{proposition}
\begin{proof} For $N\in\mathbb{N}$, let $x_N =(1/N) \sum_{i=1}^{2^N-1} h_i$ and let
$y \in \mathcal{F}_\delta((h_i))$. Note that $\|x_N\|=1$. We shall prove that $\|x-y\|\ge 1$ provided
$N\ge2/\delta$.
Since $(h_i)$ is  a monotone basis, we may assume that $y \in
\operatorname{span}\{h_i\colon 0\le i \le 2^N-1\}$.
Since $x_N$ and $-x_N$ have the same distribution,
we may also assume that  the coefficient of $h_0$ in the expansion of $y$ is $-\alpha$, where $\alpha\ge0$.
Let $k_1\ge1$
be the first level (if there are any) of the Haar system for which the \textit{leftmost} Haar function has a nonzero coefficient in the expansion of $y$. Let this Haar function be $h_{i_1}$ and let $a_1$ be the corresponding
coefficient. Note that $|a_1| \ge \delta$.
By considering the left-hand and the right-hand halves of the support
of $h_{i_1}$, and using the monotonicity of the Haar basis, we see that
$$\max_{t \in I_1}(x - y)(t) \ge \frac{k_1-2}{N} + \alpha+ \delta
= \frac{k_1}{N} +\alpha + (\delta - \frac{2}{N}),$$
where $I_1$ is  the (left-hand or right-hand)  half of the support of $h_{i_1}$ on which
$a_1h_{i_1}$ takes a \textit{negative} value.
Now we repeat the argument for $I_1$. Suppose that the next level for which
there is a nonzero coefficient in the leftmost Haar function whose support is entirely contained in $I_1$ is  the $(k_1+k_2)$-th level, where $k_2\ge1$. Let $h_{i_2}$ denote this Haar function and let $a_2$ be the corresponding coefficient. Then, by the same reasoning as above, we get
$$ \max_{t \in I_2}(x-y)(t)> \frac{k_1+k_2}{N}  +\alpha + 2(\delta - \frac{2}{N}),$$
where $I_2$  is  the  half of the support of $h_{i_2}$ on which
$a_1h_{i_2}$ takes a \textit{negative} value. This process terminates after $J\ge0$ steps
at level $k_1+\dots+k_J$  with a set $I_J$ (half of the support of $h_{i_J}$) such that
$$\max_{t \in I_J}(x-y)(t)> \frac{k_1+\dots+k_J}{N}  +\alpha + J(\delta - \frac{2}{N}).$$
Finally, let $I$ be the left-hand half of the leftmost Haar function on the $N$-th level
whose support is entirely contained in $I_J$. Since the inductive process has terminated after $J$ steps, we obtain
$$(x-y)(t)\ge 1  +\alpha + J(\delta - \frac{2}{N}) \ge 1\qquad (t \in I)$$
provided $N\ge2/\delta$.

\end{proof}
\begin{remarks} (i) The proof of Proposition~\ref{prop: haarsystem} actually shows that
if $\delta\ge2/N$ then $\mathcal{F}_\delta((h_i))$
is not an $\varepsilon$-net for the unit ball of $\ell_\infty^{2^N}$ for any $0<\varepsilon<1$. \newline
(ii) In the terminology of Section~\ref{sec: BQP} below, Proposition~\ref{prop: haarsystem} shows that $(h_i)$ does not have the Net Quantization Property.
\end{remarks}
\section{An Existence Result} \label{sec: existence}
\begin{theorem} \label{thm: existence} Suppose that $c_0\hookrightarrow X$. Then $X$ has a
bounded, total,  weakly null, \mnote{I added "weakly null" which is relevant later.} normalized
 minimal system which has the $(\varepsilon, c\varepsilon)$-SCQP for all $\varepsilon>0$,
where $c$ is an absolute constant (independent of $X$ and $\varepsilon$). Moreover, if $X$ has a basis,
then $X$ has a normalized weakly null $(\varepsilon, c\varepsilon)$-SCQP basis.
\end{theorem}
First let us
explain the construction that is used in the proof of Theorem~\ref{thm: existence}.
To that end,
let $(e_j)_{j=1}^{n+1}$ denote
the unit vector basis of $\ell_\infty^{n+1}$. Define a new basis $(f_j)_{j=1}^{n+1}$ as follows:
\begin{equation*}
f_j = e_j + \frac{e_{n+1}}{n} \qquad (1 \le j \le n)
\end{equation*}
and
\begin{equation*}
f_{n+1} = e_1+e_2+\dots+e_n. \end{equation*}
The following lemma is easily verified.
\begin{lemma} \label{lem: newbasis} $(f_j)_{j=1}^{n+1}$ is a
normalized basis for $\ell_\infty^{n+1}$ with basis constant
at most $3$. \end{lemma}
\begin{proof}[Proof of Theorem~\ref{thm: existence}] \mnote{I have added more details here.}
 By Sobczyk's theorem \cite{S} that $c_0$ is \linebreak $2$-complemented in any separable superspace and  James's theorem \cite{J} that
every Banach space  isomorphic to  $c_0$ contains an almost isometric copy of $c_0$, it follows that $X$ is
uniformly  isomorphic to $X \oplus_\infty c_0$.
So by Proposition~\ref{prop: isom} and Proposition~\ref{prop: normalize}, it suffices to prove the result for $X\oplus_\infty c_0$.
Let $(\phi_i)$ be a normalized total  minimal system (resp.\
 normalized  basis) for $X$.

For convenience, we regard $c_{00}$ as the space of all finitely supported sequences $(a^n_j)$
doubly indexed by $n \in \mathbb{N}$ and $1 \le j \le n^2+1$. Let $(e^n_j)$ denote the standard basis for this realization  of $c_{00}$ and order the basis
elements lexicographically (i.e., $e^1_1, e^1_2, e^2_1, e^2_2,\dots$).
Define a norm $\|\cdot\|_Y$ on $c_{00}$
as follows: \begin{equation*}
\|(a^n_j)\|_Y= \max\left\{\sup_{n\ge1}\|(a^n_j+a^n_{n^2+1})_{j=1}^{n^2}\|_\infty,
\left\|\sum_{n=1}^\infty \frac{1}{n^2}\left(\sum_{j=1}^{n^2}a^n_j\right)\phi_n\right\|_X\right\},
\end{equation*}  and let $Y$ denote the completion of $(c_{00},\|\cdot\|)$.
It is easily seen that
$Y$ is \textit{isometrically isomorphic} to $X\oplus_\infty (\sum_{n=1}^\infty \oplus \ell_\infty^{n^2})_0$, which in turn is
isometrically isomorphic to
$X\oplus_\infty c_0$, and that  $(e^n_j)$ is a
normalized bounded  and total minimal system for $Y$. Moreover, for each $n \in \mathbb{N}$,
 $(e^n_j)_{j=1}^{n^2+1}$ is isometrically equivalent to the basis $(f_j)_{j=1}^{n^2+1}$
described above. Thus, for  the case in which $(\phi_i)$ is a basis for
$X$, it follows easily from Lemma~\ref{lem: newbasis}  that $(e^n_j)$ is a basis for $Y$.

Let us next check that $(e^n_j)$ is weakly null. \mnote{Added proof of weakly null.}
Under the isometric isomorphism of $Y$ with $X\oplus_\infty (\sum_{n=1}^\infty \oplus \ell_\infty^{n^2})_0$,
the basis vector $e^n_j$ corresponds to \begin{equation} \label{eq: weaklynullbasis} \begin{cases}
(\phi_n/n^2,g^n_j), &\text{if $1 \le j \le n^2$}\\
(0,\sum_{i=1}^{n^2}g^n_i),  &\text{if $j = n^2+1$,} \end{cases} \end{equation} where  $(g^n_i)_{i=1}^{n^2}$ denotes the unit vector basis of $\ell_\infty^{n^2}$. Thus it sufffices to check that the sequence defined by \eqref{eq: weaklynullbasis} is weakly null. But this is readily verified directly
using the fact that
$(X\oplus_\infty (\sum_{n=1}^\infty \oplus \ell_\infty^{n^2})_{0})^*$ is isometrically isomorphic to $X^*\oplus_1 (\sum_{n=1}^\infty \oplus \ell_1^{n^2})_1$.

To see that $(e^n_j)$ has the $(\varepsilon,c\varepsilon)$-SCQP, let $\delta>0$ and
let $(D^n_j)$ be a doubly-indexed family of $\delta$-nets and
let $(a^n_j) \in c_{00}$. For each $n \in \mathbb{N}$, choose $d^n_j \in D^n_j$, with
$d^n_j=0$ if $a^n_j=0$, such that
\begin{equation} \label{eq: choiceofds1}
\left|\sum_{j=1}^k (a^n_j-d^n_j)\right| \le \delta \qquad(1 \le k \le n^2)
 \end{equation}
and
\begin{equation} \label{eq: choiceofds2}
|a^n_{n^2+1}-d^n_{n^2+1}| \le \delta. \end{equation}
>From \eqref{eq: choiceofds1} and the triangle inequality, we see that
\begin{equation} \label{eq: choiceofds3}
|a^n_j-d^n_j| \le 2 \delta \qquad (1 \le j \le n^2). \end{equation}
Combining \eqref{eq: choiceofds1},  \eqref{eq: choiceofds2}, and \eqref{eq: choiceofds3},
we obtain
\begin{equation*} \begin{split}
\|\sum (a^n_j-d^n_j)e^n_j\|_Y &\le
\sup_n\max_{1\le j \le n^2} |a^n_j-d^n_j+a^n_{n^2+1}-d^n_{n^2+1}|\\
&\quad +\sum_{n=1}^\infty\frac{1}{n^2} |\sum_{j=1}^{n^2}(a^n_j-d^n_j)|\\
&\le 3\delta  +\delta\cdot \frac{\pi^2}{6}. \end{split}
\end{equation*}
This shows that $(e^n_j)$ is a minimal system (resp.\ basis) for $Y$ with the $(\varepsilon,c\varepsilon)$-SCQP for $c=(3+\pi^2/6)^{-1}$.
\end{proof}
\begin{remark} The construction used in the proof of Theorem~\ref{thm: existence}
was first used by Wojtaszczyk \cite{W2}. The dual construction was
used recently in \cite{DKK} to construct a \textit{quasi-greedy}
basis for $L_1([0,1])$.\end{remark}

\section{The Net Quantization Property}\label{sec: BQP}
In this section we discuss a  natural quantization property which is more general than the CQP.
\begin{definition} \label{def: BQP} \mnote{This is  basically the old BQP but now for the whole space}
Let $\varepsilon>0$ and let $\delta>0$. \newline
(a) A dictionary $(e_i)$ has the $(\varepsilon,\delta)$-\textit{Net Quantization Property}
(abbr.\ $(\varepsilon,\delta)$-NQP) if  for every $x\in X$ there exist a finite subset
$E \subset \mathbb{N}$
and $n_i \in \mathbb{Z}$ ($i \in E$) such that  \begin{equation} \label{eq: NQPdef}
\|x - \sum_{i\in E} n_i\delta e_i\| \le \varepsilon. \end{equation}
(b)  $(e_i)$ has the NQP if $(e_i)$ has the
$(\varepsilon,\delta)$-NQP for some $\varepsilon>0$ and $\delta>0$.
\end{definition}

\begin{remarks} (i)  Note that \eqref{eq: NQPdef}  simply says that $\mathcal F_\delta((e_i))$ is an $\varepsilon$-net for $X$.
In particular,
choosing $x= \sum_{i \in F} a_ie_i$ in \eqref{eq: NQPdef},  it is important to emphasize that
the set $E$ is \textit{not} required to be contained in $F$. This suggests that the
the NQP  property should be weaker than the CQP property, and we prove below that this is indeed the case.   \newline
(ii) The analogue of Proposition~\ref{prop: elemproperties} remains valid for the NQP.
\end{remarks}
The analogue of Theorem~\ref{thm: ballenough} for the NQP which is stated below remains valid with essentially the same proof.
\begin{theorem} \label{thm: ballenoughBQP} Let $(e_i)$ be a dictionary for $X$. The following are equivalent:
\begin{itemize} \item[(i)] $(e_i)$ has the NQP;
\item[(ii)] there exist $0<\varepsilon<1$ and $\delta>0$ such that $\mathcal F_\delta((e_i))$ is an $\varepsilon$-net for $Ba(X)$.
\end{itemize} \end{theorem}
\begin{corollary} Let $X$ be a separable Banach space. There exists a  dictionary $(e_i)$ with the NQP such that $\mathcal{F}_1((e_i))$ is
$M$-dense in $X$ and $(1/M)$-separated for some $M>0$.
 \end{corollary}
\begin{proof} Let $(x_n)_{n=1}^\infty$ be a semi-normalized fundamental bounded minimal system for $X$ with $\|x_n\| \le 1/3$ for all $n$. Let $(y_n)$ be  dense  in the unit ball of $X$  with $y_n \in \langle x_i \rangle_{i=1}^{n-1}$, and let $e_n = x_n+y_n$. Then $(e_n)$ is semi-normalized and $1/2$-dense in $Ba(X)$. So by Theorem~\ref{thm: ballenoughBQP}
$\mathcal{F}_1((e_i))$ is an $M$-net for $X$ for some $M>0$. Using the fact that $(x_i)$ is a bounded minimal system it is easily verified that
 $\mathcal{F}_1((e_i))$ is  $(1/M)$-separated for sufficiently large $M$.
\end{proof}
The counterpart to Corollary~\ref{cor: preciseversion} takes the following form. This result seems to be of interest even when $X$ is finite-dimensional.
\begin{theorem}\label{cor: preciseversion2} Let $0<\varepsilon_0<1$, $\delta>0$, and let $(e_i)$ be a (not necessarily semi-normalized) fundamental system for $X$. If $\mathcal{F}_{\delta}((e_i))$ is $\varepsilon_0$-dense in $Ba(X)$
then   $\mathcal{F}_{\delta}((e_i))$ is $\varepsilon_1$-dense in $X$ for all
$$\varepsilon_1> (\left\lfloor \frac{\varepsilon_0}{1-\varepsilon_0}\right\rfloor+1)\varepsilon_0.$$
In particular, if $\varepsilon_0<1/2$, then $\mathcal{F}_{\delta}((e_i))$ is $\varepsilon_1$-dense in $X$ for all $\varepsilon_1>\varepsilon_0$.
\end{theorem}

Next we introduce the analogue of the SCQP.
\begin{definition} \label{def: SNQP} \label{This is the "strong" version}
Let $\varepsilon>0$ and let $\delta>0$. \newline
(a) A dictionary $(e_i)$ has the $(\varepsilon,\delta)$-\textit{Strong Net Quantization Property}
(abbr.\ $(\varepsilon,\delta)$-SNQP) if $\mathcal F_{\overline{D}}((e_i))$ is an $\varepsilon$-net for $X$
 for every sequence $\overline{D} = (D_i)$ of $\delta$-nets.  \newline
(b)  $(e_i)$ has the SNQP if $(e_i)$ has the
$(\varepsilon,\delta)$-SNQP for some $\varepsilon>0$ and $\delta>0$.
\end{definition}
The proof of Proposition~\ref{prop: uniformity} for a general dictionary does not seem to transfer to the SNQP. However,
when $(e_i)$ is a Schauder basis it is easy to modify the proof
 to get the following uniform boundedness result.
\begin{proposition} Let $(e_i)$ be a Schauder basis for $X$. The following are  equivalent:
\begin{itemize}
\item[(i)] $(e_i)$ has the SNQP;
\item[(ii)] for all $\delta>0$ and for every sequence  $\overline{D}=(D_i)$ of $\delta$-nets there exists $M := M(\overline{D})>0$ such that
$\mathcal F_{\overline{D}}((e_i))$ is an $M$-net for $X$;
 \item[(iii)] condition (ii) for $\delta=1$.
\end{itemize} \end{proposition}
\begin{remark} The analogues of Propositions~\ref{prop: isom} and \ref{prop: normalize}  remain valid
for the SNQP. \end{remark}

Trivially, every separable Banach space  has a dictionary with
the $(\varepsilon,c\varepsilon)$-SNQP for all $0<c<1$. Indeed, simply
take
$(e_i)$ to be dense in the unit sphere of $X$. By
a more careful choice of  dense set in the unit sphere of $\ell_2$,  it is not difficult
to construct an  NQP dictionary for $\ell_2$ which is \textit{not} a CQP dictionary.
Our next result, the construction of  an SNQP {\it Schauder basis} which  does not have the CQP, is more involved.
It is a consequence of the following general embedding theorem. (Recall that a Schauder basis is \textit{bimonotone}
if the basis projections $(P_n)$ satisfy $\|P_n\| = \|I-P_n\| =1$ for all $n \ge 1$.)
\begin{theorem} \label{thm: NQPnotCQP} \mnote{This is Ted's embedding theorem essentially}
Let $(e_i)$ be a normalized bimonotone basis for a Banach space $E$. Given $\eta>0$ there exists a Banach space $U$
with a normalized  monotone basis $(u_i)$ with the following properties: \begin{itemize}
\item[(a)] $(u_i)$ has the $(\varepsilon, \varepsilon/3)$-SNQP; \\
\item[(b)] there exists a subsequence $(u_{n_i})$ of $(u_i)$ that is $(1+\eta)$-equivalent to $(e_i)$.
\end{itemize} \end{theorem}
Before proceeding with the proof, let us see how it implies the existence of an SNQP basis which is not a CQP basis.
The CQP is inherited by subsequences, so if we apply Theorem~\ref{thm: NQPnotCQP} to any basis $(e_i)$ which does not have the CQP
(e.g. the unit vector basis of $\ell_2$) then the constructed basis $(u_i)$ will have the SNQP  but not the CQP.

\begin{proof}[Proof of Theorem~\ref{thm: NQPnotCQP}] Choose integer reciprocals $\eta_i  \downarrow 0$  such that for each $j$ the set
 $$S_j:=\{\sum_{i=1}^j k_i \eta_i e^*_i \colon k_i \in \mathbb{Z}\}\cap Ba(E^*)$$ is $(1-\eta)$-norming for $\langle e_i \rangle_{i=1}^j$.
Note that if $j \le k$, then each element $g$ of $S_k$ is an ``extension'' of an element $g'$ of $S_j$ (i.e. $g(e_i) = g'(e_i)$ for $1 \le i \le j$).
Note also that $(e^*_i)_{i=1}^j \subset S_j$ since $(e_i)$ is bimonotone.
We shall construct  a subset $\mathcal G \subset Ba(c_{0})\cap c_{00}$  such that
$P_n(\mathcal G) \subset \mathcal G$ for all $n \in \mathbb{N}$,
where $(P_n)$ is the sequence of basis projections in $c_{00}$. Then we define $U$  to be  the Banach space with Schauder basis
$(u_i)$ whose norm is given by
$$ \|\sum a_i u_i\| = \sup_{f \in \mathcal G} |\sum f(i) a_i|.$$
The conditions on $\mathcal G$ ensure that $(u_i)$ is a monotone basis for $U$. The construction of $\mathcal G$ and the sequence $(n_i)$ is inductive.
Set $n_1 = 1$ and
$$\mathcal G_1:= \{ (k_1 \eta_1,0,0,\dots)  \colon k_1 \eta_1 e_1^* \in S_1\}.$$ Suppose $j_0 \ge 1$ and that $n_j$ and $\mathcal G_j$ have been defined for each $j \le j_0$ such that every $f \in \mathcal G_j$ is supported on $[1,n_j]$, $P_n(\mathcal G_j) \subset \mathcal G_j$ for all $n \in \mathbb{N}$, and
  $P_{n_j}(\mathcal G_{j+1}) \subset \mathcal G_j$, i.e. every element of $\mathcal G_{j+1} \setminus \mathcal G_j$ is an extension on $[n_j+1, n_{j+1}]$ of
 some element of $\mathcal G_j$, and such that if $f \in \mathcal G_j$ then there exists a  $\tilde g := \tilde g(f)  \in S_j$
such that $f(n_i) = \tilde g(e_i)$ for all $1 \le i \le j$ (and, conversely, for every $g \in S_j$ there exists
$f \in \mathcal G_j$ such that $g = \tilde g(f)$).
We now proceed to the definition of $n_{j_0+1}$ and $\mathcal G_{j_0+1}$.
Let $$T_{j_0} := \{(f,g) \in \mathcal G_{j_0} \times S_{j_0+1} \colon \text{ $g$ extends $\tilde g(f)$} \} \subset \mathcal G_{j_0} \times S_{j_0+1}.$$
Let $n_{j_0+1} := n_{j_0} + \operatorname{card}{T_{j_0}} + 1$ and define a bijection $(f,g) \rightarrow i((f,g))$ from $T_{j_0}$ onto $[n_{j_0}+1, n_{j_0+1}-1]$.
For each $(f,g) \in T_{j_0}$, define $f' := f'((f,g))$ by \begin{equation*}
f'(i) = \begin{cases} f(i) &\text{if $1 \le i \le n_{j_0}$}\\
                      1     &\text{if $i= i((f,g))$}\\
                      g(e_{j_0+1}) &\text{if $i = n_{j_0+1}$}\\
                      0 &\text{otherwise.} \end{cases} \end{equation*}
Set $$\mathcal G_{j_0+1} := \{P_n(f'((f,g))) \colon (f,g) \in T_{j_0}, n \le n_{j_0+1} \}.$$ Finally, define
 $\mathcal G=\cup_{j \ge 1} \mathcal G_j$.
Then  $\mathcal G$  satisfies
$P_n(\mathcal G) \subset \mathcal G$ ($n \in \mathbb{N}$) as claimed. Thus,
$(u_i)$ is a  monotone basis for $U$. Moreover, since $e_j^* \in \mathcal{S}_j$, it is easily checked that $\|u_i\|=1$ for all $i$.
Henceforth, we identify $\mathcal G$ with a norming subset of $Ba(U^*)$ and use the notation
$f(\sum a_i u_i) := \sum f(i)a_i$ for $f \in \mathcal G$.
It is clear from the construction that
$$\|\sum_{i=1}^{m} a_i u_{n_i}\| = \sup_{g \in S_m} g(\sum_{i=1}^m a_i e_i),$$
and so $(u_{n_i})$ is $(1+\eta)$-equivalent to $(e_i)$, which verifies (b).

Let us now turn to the verification of (a). Let $\varepsilon>0$ and let $(D_i)$ be a sequence of
$\varepsilon/3$-nets. To show that $\mathcal F_{\overline{D}}((u_i))$ is an $\varepsilon$-net for $U$, it suffices
to show that for every
 $ x = \sum_{i \in A} a_i u_i \in U$, where $A\subset \mathbb{N}$ is finite, there exists
$y = \sum_{i \in E} d_i e_i$ ($d_i \in D_i$), where  $E \subset \mathbb{N}$ is finite,
 such that  $\|x - y\| \le 2\varepsilon/3$ (since the collection of all such $x$ is dense in $U$). We may assume that $A \subset [1, n_j]$ for some $j$. The proof is by induction on $j$. The  case $j=1$ is clear: $n_1=1$, so $x = a_1u_1$ in this case, and we simply choose
$d_1 \in D_1$ with $|a_1-d_1| \le \varepsilon/3$, so that $\|x-d_1u_1\| \le \varepsilon/3$.

Suppose the inductive hypothesis holds for $j = j_0$. For the  inductive step, suppose that  $x = \sum_{i=1}^{n_{j_0+1}} a_i u_i$ and let $x' = \sum_{i=1}^{n_{j_0}} a_iu_i$. By the inductive hypothesis there exists
$y' = \sum_{i=1}^{n_{j_0}} d_i u_i$ such that $\|x'-y'\| \le 2\varepsilon/3$. Let $y = \sum_{i=1}^{n_{j_0+1}}d_iu_i$ be an extension of $y'$ to $[1,n_{j_0+1}]$.
Then
$$ |f(x-y)| = | f(x'-y')| \le 2\varepsilon/3 \quad \text{for all $f \in \mathcal G_{j}$ when $j \le j_0$}.$$
Since $P_{n_{j_0+1}}(\mathcal G_j) = \mathcal G_{j_0+1}$ when $j \ge j_0+1$, it suffices to choose the extension $y$ such that $|f(x-y)| \le 2\varepsilon/3$
for all $f \in \mathcal G_{j_0+1} \setminus \mathcal G_{j_0}$. To that end, for each $(f,g) \in T_{j_0}$,
setting $i' := i((f,g))$ choose
 $d_{i'} \in D_{i'}$ such that $$|f(x'-y') +a_{i'} - d_{i'}| \le \varepsilon/3.$$ This
defines
 $d_i$  for $n_{j_0}+1 \le i \le n_{j_0+1} -1$. Finally,
choose $d_{n_{j_0+1}}\in D_{n_{j_0+1}}$ such that $|a_{n_{j_0+1}}-d_{n_{j_0+1}}| \le \varepsilon/3$. This completes the definition of $y$.
Suppose that $f'=f'((f,g))$ for some $(f,g)\in T_{j_0}$.
Then \begin{align*}
|f'(x-y)| &= |f(x'-y') + a_{i'} - d_{i'} + g(e_{j_0+1})(a_{n_{j_0+1}}-d_{n_{j_0+1}})| \\
&\le |f(x'-y') + a_{i'} - d_{i'}| + |(a_{n_{j_0+1}}-d_{n_{j_0+1}})|\\ &\le \varepsilon/3 + \varepsilon/3 = 2\varepsilon/3.
\end{align*}
Moreover,
$$|(P_nf')(x-y)| =|f(x'-y') + a_{i'} - d_{i'}| \le \varepsilon/3 \qquad(i'\le n \le n_{j_0+1}-1)$$
and
$$|(P_nf')(x-y)| =|(P_nf)(x'-y') | \le 2\varepsilon/3 \qquad(1 \le n < i')$$
 by inductive hypothesis since $P_nf \in \mathcal{G}_{j_0}$ for $1 \le n < i'$.
This completes the proof of the inductive step.
\end{proof}
In some of our results in Section~\ref{sec: containmentofc0} it is possible to replace the CQP  by the formally weaker assumption that every subsequence has the NQP. When $(e_i)$ is a Schauder basis, however, our next result shows  that this assumption is in fact equivalent to the CQP.
\begin{theorem} Let $(e_i)$ be a semi-normalized basic sequence which fails the CQP. Then some subsequence fails the NQP for its closed linear span. \end{theorem}
\begin{proof} Let $K$ be the basis constant of $(e_i)$. We may assume without loss of generality that $\|e_i\| \le 1$ for all $i$.

\textit{Claim 1:} For every $\delta>0$ there exists $M \subset \mathbb{N}$ such that $(e_i)_{i \in M}$ fails the $(1, \delta)$-NQP.

\textit{Proof of Claim~1:} Suppose not. Then there exists $\delta>0$ such that $(e_i)_{i \in M}$ has the $(1, \delta)$-NQP
for every $M \subset \mathbb{N}$. Let $x = \sum_{i \in E} a_i e_i$ and let $n = \max E$. Since $M := E \cup (n,\infty)$ has the $(1, \delta)$-NQP
there exists $y \in \mathcal{F}_\delta((e_i)_{i \in M})$ such that $\|x-y\| \le 1$. Then $\|x - P_E y\|  \le K$. Thus $(e_i)$ has the
$(K,\delta)$-CQP, which is a contradiction.

\textit{Claim 2:} For all $n \in \mathbb{N}$ there exist a finite set $F_n \subset [n+1, \infty)$ and $x_n = \sum_{i \in F_n} a_i e_i$
such that $\|y - x_n\| > 2K$ for all $y \in \mathcal{F}_{1/n}((e_i)_{i \in F_n})$.

\textit{Proof of Claim~2:} Let $\delta_n = 1/n$. By Claim~1 there exists $M_n \subset \mathbb{N}$ such that $(e_i)_{i \in M_n}$ fails the
$(2K+1, \delta_n)$-NQP. So there exists $z_n = \sum_{i \in M_n} a_i e_i$ with $\|z_n - y\| > 2K+1$ for all
$y \in \mathcal{F}_{1/n}((e_i)_{i \in M_n})$. Let $x_n = z_n|_{[n+1,\infty)}$. Note that every vector supported on $[1,n]\cap M_n$ (in particular, the vector
$z_n-x_n$)
can be $1$-approximated by an element of $\mathcal{F}_{1/n}((e_i)_{i \in M_n})$ simply by approximating each of the (at most $n$) nonzero coordinates to within $\delta_n = 1/n$. Setting $F_n := \supp x_n$, it
 follows that $\|x_n-y\|>2K$ for all $y \in \mathcal{F}_{1/n}((e_i)_{i \in F_n})$. Thus, $x_n$ and  $F_n$ verify Claim~2.

Now pass to a subsequence so that the sets $F_{n_k}$ satisfy $\max F_{n_k} < \min F_{n_{k+1}}$ for all $k \in \mathbb{N}$. Let $M = \cup_{k \ge 1}
F_{n_k}$.

\textit{Claim~3:} $(e_i)_{i \in M}$ fails the NQP.

\textit{Proof of Claim~3:} Suppose that $(e_i)_{i \in M}$ has the $(1,\delta)$-NQP (and hence the $(1, 1/n)$-NQP provided $1/n < \delta$).
Choose $k$ with $1/n_k < \delta$. Then there exists $y \in \mathcal{F}_{1/n_k}((e_i)_{i \in M})$ such that $\|y - x_{n_k}\| \le 1$.
But this implies that $\|P_{F_{n_k}}(y) - x_{n_k}\| \le 2K$, which contradicts the choice of $x_{n_k}$ and $F_{n_k}$.

\end{proof}
We turn now to discuss the relationship between the NQP and unconditionality.

\begin{theorem} \label{thm: NQPnotunc}
Suppose that $X$ has a semi-normalized unconditional basis $(e_i)$ with the NQP. Then $(e_i)$ is equivalent to the unit vector basis of
$c_0$. \end{theorem}
\begin{proof} Let $K$ be the constant of unconditionality of $(e_i)$ and choose $\varepsilon >0$ such that $K< \dfrac{1-\varepsilon}{\varepsilon}$.
There exists $\delta>0$ such that $\mathcal F_\delta((e_i))$ is $\varepsilon$-dense in $X$. Suppose $x = \sum e^*_i(x) e_i \in X$ with $\|x\| = 1$ and
$\|x\|_\infty := \sup |e^*_i(x)| = \alpha<\delta$. Choose $y \in \mathcal F_\delta((e_i))$ with $\|x - y\| \le \varepsilon$. Then
 $$\|y\| \ge \|x\|-\|x-y\| \ge 1 - \varepsilon.$$
Since $\sup |e^*_i(x)| \le \alpha$ and since $y \in \mathcal F_\delta((e_i))$, it follows  that $y = \sum \lambda_i e^*_i(x-y)e_i$ for a multiplier sequence $(\lambda_i)$
satisfying $$\sup |\lambda_i| \le \frac{\delta}{\delta-\alpha}.$$
 Hence by $K$-unconditionality of $(e_i)$, we have \begin{equation} \label{eq: Kcontra}
1-\varepsilon \le \|y\| \le K \sup |\lambda_i| \|y-x\| \le K \frac{\delta}{\delta-\alpha}\varepsilon. \end{equation}
If $(e_i)$ is not equivalent to the unit vector basis of $c_0$ then  $\alpha$ may be chosen to be arbitrarily small.
But then \eqref{eq: Kcontra} yields $K \ge \dfrac{1-\varepsilon}{\varepsilon}$, which contradicts the choice of $\varepsilon$.

\end{proof}

Weaker notions of unconditionality (see \cite{DOSZ}), especially that of a {\it quasi-greedy} basis, have recently attracted attention in connection
 with greedy algorithms for data compression.
  Our next goal is to show that every  quasi-greedy basis with the NQP is equivalent to the unit vector basis of $c_0$.
The relevant definitions are given next. For further information on the topic of  greedy algorithms in Banach spaces, we refer the reader to
\cite{KT, DKK, DKKT, DOSZ, W1}. \mnote{The rest of this section has been rewritten to improve "quasi-greedy" to "Elton-unconditional" which is more general. The proofs are different.}

\begin{definition} Let $(e_i)$ be a dictionary for $X$ and let $\delta>0$.
\newline
(a) Denote by $L((e_i), \delta)$ the least constant
$L \in [1,\infty]$ with the property that whenever $\|\sum a_ie_i\| \le 1$ and
$F \subset \{i\colon |a_i| \ge \delta\}$ then
\begin{equation*}
\|\sum_{i \in F} a_i e_i\| \le L.
\end{equation*}
(b) We say that $(e_i)$ is \textit{Elton-unconditional} if
$$L((e_i), \delta)< \infty \quad \text{for all
$\delta>0$}.$$
\newline
(c) Denote by $K((e_i), \delta)$ the least constant
$K \in [1,\infty]$ with the property that whenever $\|\sum a_ie_i\| \le 1$ and
$F = \{i\colon |a_i| \ge \delta\}$ then
\begin{equation*}
\|\sum_{i \in F} a_i e_i\| \le K.
\end{equation*}
(d) We say that $(e_i)$ is \textit{quasi-greedy} if \begin{equation*}
K((e_i)) := \sup_{\delta>0} K((e_i), \delta) < \infty. \end{equation*}
\end{definition}
\begin{remark} Clearly, $K((e_i),\delta) \le L((e_i),\delta))$. Note that $(e_i)$ is unconditional if and only
if $\sup_{\delta>0} L((e_i), \delta)< \infty$. It is known that every quasi-greedy basic  sequence is Elton-unconditional
(in fact a semi-normalized Schauder basis $(e_i)$ is Elton-unconditional if and only if $K((e_i), \delta)< \infty$ for all $\delta>0$)
and that there exist Elton-unconditional bases which are not quasi-greedy \cite{DKK}.
\end{remark}
\begin{lemma} Let $(e_i)$ be a minimal system for $X$. Suppose that there exist $0< \varepsilon <1$, $\delta>0$, and $\lambda>0$ such that
$\mathcal F_\delta((e_i))\cap \lambda Ba(X)$ is an $\varepsilon$-net for $Ba(X)$ and such that $L((e_i), \delta/\lambda)< \infty$.
Then $(e_i)$ is equivalent to the unit vector basis of $c_0$.
\end{lemma}
\begin{proof} Clearly, $(e_i)$ has the NQP. So by Theorem~\ref{thm: NQPnotunc} it suffices to show that $(e_i)$ is unconditional.
Let $S := \mathcal F_\delta((e_i)) \cap \lambda Ba(X)$. Since $S$ is an $\varepsilon$-net for $Ba(X)$ it follows
that $S$ is
$(1-\varepsilon)$-norming for $X^*$, i.e.
$$\|x^*\| \le \frac{1}{1-\varepsilon} \sup \{|x^*(x)| \colon x \in S\} \qquad(x^* \in X^*).$$
Moreover, if $x=\sum_E k_i\delta e_i \in S$ and $F \subseteq E$, then (since  $x/\lambda \in Ba(X)$)
$\|\sum_{i \in F} k_i \delta e_i\| \le \lambda L((e_i), \delta/\lambda)$. Hence
$$\tilde S := \{ \sum_{i \in F} k_i\delta e_i \colon \sum_{i \in E} k_i\delta e_i \in S, F \subseteq E\}  \subseteq \lambda L((e_i), \delta/\lambda) Ba(X).$$
Now suppose that $\sum_{i \in E} a_ie_i^* \in X^*$ and that $F \subseteq E$. Then
\begin{align*} \|\sum_{i \in F} a_ie_i^*\| &\le \frac{1}{1-\varepsilon} \sup \{ \sum_{i \in F} a_ie^*_i(x) \colon x \in S\}\\
&\le \frac{1}{1-\varepsilon} \sup \{ \sum_{i \in E} a_ie^*_i(x) \colon x \in \tilde S\}\\
& \le \frac{\lambda}{1-\varepsilon} L((e_i), \delta/\lambda) \|\sum_{i \in E} a_i e_i^*\|. \end{align*}
Thus, $(e^*_i)$ is $K$-unconditional for $K= \lambda L((e_i), \delta/\lambda)/(1-\varepsilon)$, and hence by duality $(e_i)$ is also $K$-unconditional.
\end{proof}
The following substantial strengthening of Theorem~\ref{thm: NQPnotunc} is an immediate consequence of the last result.
\begin{theorem} Suppose that $(e_i)$ is a minimal system with the NQP. If $(e_i)$ is Elton-unconditional (in particular, if $(e_i)$ is quasi-greedy) then
$(e_i)$ is equivalent to the unit vector basis of $c_0$. \end{theorem}

The main open question of this section is the following.

\begin{problem} \label{prob: NQPcontainc0} Suppose that $X$ has an NQP basis. Does $c_0 \hookrightarrow X$? \end{problem}

In fact, we do not know whether or not $\ell_1$ provides a negative answer to Problem~\ref{prob: NQPcontainc0}.
\begin{problem} \mnote{It would be nice to know the answer here.}
Does $\ell_1$ have an NQP basis (resp.\ minimal system)? \end{problem}
 We conclude this section with some partial results concerning Problem~\ref{prob: NQPcontainc0}.
\begin{theorem} \label{thm: embedl1indual}
Suppose that $(e_i)$ is a bounded NQP minimal system for $X$.  Then $\ell_1 \hookrightarrow X^*$. In fact, either $\ell_\infty \hookrightarrow X^*$
or a subsequence of $(e_i^*)$
is equivalent to the unit vector basis of $\ell_1$. \end{theorem}
\begin{proof} Since $(e_i)$ has the NQP,  there exists $\delta>0$
such that $\mathcal{F}_\delta((e_i))\cap (3/2)Ba(X)$ is a $1/2$-net for $Ba(X)$.
Thus, \begin{equation*}
\frac{1}{2}\|x^*\| \le \sup\{|x^*(x)| \colon x \in \mathcal{F}_\delta((e_i))\cap \frac{3}{2}Ba(X)\}
\le \frac{3}{2}\|x^*\|\qquad(x^* \in X^*). \end{equation*}
So $X^* \hookrightarrow C(K)$, where $K$ is the weak-star closure of \linebreak $\mathcal{F}_\delta((e_i)) \cap (3/2) Ba(X)$ in $X^{**}$.
By Rosenthal's $\ell_1$ theorem \cite{R1}, $(e^*_i)$ has a subsequence equivalent to the unit vector basis of $\ell_1$ or a weakly Cauchy subsequence $(f_i^*)$. The former
obviously implies that $\ell_1 \hookrightarrow X^*$.
In the latter case, let
$g_i^* = f_{2i}^*-f^*_{2i-1}$. Then $(g_i^*)$ is weakly null in $X^*$, and, since the range of $f_i^*|_K \subset \mathbb{Z}\delta$, we have
 \begin{equation*}
g_i^*(k) \ne0 \Rightarrow |g^*_i(k)| \ge \delta \qquad(i \in \mathbb{N}, k \in K).
\end{equation*}
Thus, for each $k \in K$, the sequence $(g^*_i(k))$ is eventually zero,  so
  the series $\sum g^*_i$ is \textit{extremely weakly unconditionally Cauchy},
i.e.\ $\sum_{i=1}^\infty |g^*_i(k)|$  converges (trivially!) for every $k \in K$.
By a theorem of Elton \cite{E2} (see also \cite{HOR}),
$c_0 \hookrightarrow [(g^*_i)]$. But this implies that $\ell_\infty$ is isomorphic to a subspace of $X^*$ \cite{BP1}, and a fortiori that $\ell_1 \hookrightarrow X^*$.
 \end{proof}
\begin{corollary} Suppose that $X$ is reflexive. Then $X$ does not contain a
bounded minimal system
with the NQP. \end{corollary}
\section {Containment of $c_0$} \label{sec: containmentofc0} \mnote{This section is new.}
The main result of this section is the following converse to Theorem~\ref{thm: existence}.
\begin{theorem} \label{thm: conversemainthm} \mnote{This is the main result.}
Let $(e_i)$ be a semi-normalized basic sequence with the CQP.
Then $(e_i)$ has a subsequence that is equivalent to the unit vector basis of $c_0$ or to the
summing basis of $c_0$. \end{theorem}
As the proof is quite long we shall break it down into  several parts. We shall  frequently refer to the
excellent survey article  \cite{AGR} for the
proofs of certain assertions.

 First we prove a
result  which is of independent interest.
\begin{thm}  \label{thm: indepinterest} \mnote{This is one of Ted's faxes}
 Let $(e_i)$ be a semi-normalized  nontrivial weakly Cauchy basis for $X$. Then there exists
a subsequence $(e_{n_i})$ such that either \begin{itemize}
\item[(a)] $(e_{n_i})$ is equivalent to the summing basis of $c_0$, or \item[(b)] $(e_{n_i}^*)$ is weakly null in $[(e_{n_i})]^*$.
\end{itemize}
\end{thm}
\begin{proof} Let $x^{**} \in X^{**} \setminus X$ be the weak-star limit of $(e_i)$.
By passing to a subsequence we may assume that $(e_i)$ dominates the summing basis, i.e.
$\|\sum a_i e_i\| \ge c\|(a_i)\|_{sb}$ for some $c>0$ \cite[Prop. II.1.5]{AGR}.
 If $x^{**} \in X^{**} \setminus D(X)$, where $D(X)$ denotes the collection of all elements of $X^{**}\setminus X$
 whose restrictions to $Ba(X^*)$ (equipped with the weak-star topology)
are differences of semi-continuous functions (see \cite{AGR}), then by
 \cite[Theorem~1.8]{R2} $(e_i)$ has a \textit{strongly summing} subsequence $(e_{n_i})$. In particular,  $(\sum_{i=1}^m e^*_{n_i})$
is a nontrivial weakly Cauchy sequence in $[(e_{n_i})]^*$ \cite[Lemma II.2.6]{AGR}, and so $(e^*_{n_i})$ is weakly null in $[(e_{n_i})]^*$, which yields (b).

Now suppose that
$x^{**} \in D(X)$. Then there exists a sequence $(x_i) \subset X$ that is equivalent to the summing basis of $c_0$  such that $x_i \rightarrow x^{**}$ weak-star
\cite[Theorem II.1.2]{AGR}.
Note that $(e_i-x_i)$ is weakly null.
If some subsequence of $(e_i - x_i)$ is norm-null then (a) follows by a standard perturbation argument.
 So we may assume that $(e_i - x_i)$ is a semi-normalized
weakly null sequence. By a theorem of Elton \cite{E1,O}, $(e_i - x_i)$ has either a subsequence equivalent to the unit vector basis of $c_0$ or
a basic subsequence whose sequence of biorthogonal functionals is weakly null (in the dual of the closed linear span of that basic subsequence). If the first alternative holds, let $(e_{n_i}-x_{n_i})$ be the $c_0$ subsequence. Then \begin{align*}
c\|(a_i)\|_{sb} &\le \|\sum a_i e_{n_i}\|\\ &\le \|\sum a_i (e_{n_i}-x_{n_i})\| + \|\sum a_i x_{n_i}\|\\ &\le C_1 \sup_i|a_i| + C_2\|(a_i)\|_{sb}
\le C_3 \|(a_i)\|_{sb}, \end{align*}
for certain constants $C_1$, $C_2$, $C_3$.
Hence $(e_{n_i})$ is equivalent to the summing basis of $c_0$. If the second alternative holds, let $(e_{n_i}-x_{n_i})$ be a basic subsequence with weakly null biorthogonal functionals.
To prove that $(e^*_{n_i})$ is weakly null in $[(e_{n_i})]^*$, it suffices to show that $a_i \rightarrow 0$ whenever $(a_i)$ satisfies
$\sup_m \|\sum_{i=1}^m a_ie_{n_i}\| =K < \infty$.
 Now
$$ \|\sum_{i=1}^m a_i x_{n_i}\| \le C_2 \|(a_i)\|_{sb}  \le c^{-1}C_2 \|\sum_{i=1}^m a_ie_{n_i}\| \le c^{-1}C_2K,$$
and hence by the triangle inequality
$$\sup_m \|\sum_{i=1}^m a_i(e_{n_i}-x_{n_i})\| \le K+c^{-1}C_2K.$$ Since the sequence of biorthogonal functionals to $(e_{n_i}-x_{n_i})$ is weakly null, we deduce finally that
$a_i \rightarrow 0$.
\end{proof}

\begin{proposition} \label{prop: nosubseqweaklynull} Suppose $X$ has a minimal system $(e_i)$ with the NQP. Then no subsequence of $(e_i^*)$ is weakly null.
\end{proposition}
\begin{proof} Let $0<\varepsilon<1$. There exists $\delta>0$
such that $\mathcal{F}_\delta((e_i))\cap 2Ba(X)$ is an $\varepsilon$-net for $Ba(X)$.
Thus, \begin{equation*}
\frac{1-\varepsilon}{2}\|x^*\| \le \sup\{|x^*(x)| \colon x \in \mathcal{F}_\delta((e_i))\cap 2Ba(X)\}
\le 2\|x^*\|\qquad(x^* \in X^*). \end{equation*}
So $X^* \hookrightarrow C(K)$, where $K$ is the weak-star closure of $\mathcal{F}_\delta((e_i))\cap 2Ba(X)$ in $X^{**}$.
 Suppose that $(e_{n_i}^*)$ is a weakly null subsequence of
$(e_i^*)$, whence $\sup_i \|e^*_{n_i}\| = C<\infty$. Thus, $$\{|e^*_{n_i}(k)| \colon k \in K\} \subset \{0\} \cup [\delta, 2C]\qquad(i\ge1),$$
so $(e_{n_i}^*)$ has an unconditional basic subsequence \cite[Theorem~23]{DOSZ} (see also \cite{GOW} and \cite{LAT}). Relabel this unconditional subsequence as $(e^*_{n_i})$ and let
$Y := [(e^*_{n_i})] \subset X^*$. Observe that $(e_{n_i}|_Y)$ is a semi-normalized unconditional basic sequence in $Y^*$ whose biorthogonal sequence
is  $(e_{n_i}^*) \subset Y$.
We claim that $(e_{n_i})$ has the NQP for its closed linear span in $Y^*$.  To prove the claim, let $x = \sum_{i \in A} a_i e_{i}$, where $A \subset \{n_i \colon i\ge1\}$
is finite. Since $(e_i)$ has the NQP for $X$ there exists $y = \sum_{i \in B} m_i \delta e_i$ with $\|x-y\| \le \varepsilon$, where $B \subseteq \mathbb N$ is finite and $m_i \in \mathbb{Z}$ for each $i$.
Let $z = \sum_{i \in B'} m_i e_i$, where $B' = B \cap \{n_i \colon i \ge1\}$. Then $y|_Y = z|_Y$ and $$\|x-z\|_{Y^*} \le \|x-y\| \le \varepsilon,$$ which proves the claim. Since $(e_{n_i}|_Y) \subset Y^*$ is an unconditional basic sequence with the NQP, it follows
from Theorem~\ref{thm: NQPnotunc} that $(e_{n_i}|_Y)$ is equivalent to the unit
 vector basis of $c_0$. But this implies that $(e^*_{n_i})$ is equivalent to the unit vector basis of $\ell_1$, which contradicts the assumption that
$(e^*_{n_i})$ is weakly null! \end{proof}

\begin{proposition} \label{prop: c0subseq} Suppose that $(e_i)$ is a weakly null dictionary for $X$. If every subsequence of $(e_i)$ has the NQP for its closed linear span
(in particular, if $(e_i)$ has the CQP) then
$(e_i)$ has a subsequence equivalent to the unit vector basis of $c_0$.
\end{proposition}
\begin{proof} By the aforementioned theorem of  Elton $(e_i)$ has a subsequence equivalent to the unit vector basis of $c_0$ or a basic subsequence $(e_{n_i})$ such that
$(e^*_{n_i})$ is weakly null in $[(e_{n_i})]^*$. But the latter cannot happen by Proposition~\ref{prop: nosubseqweaklynull}.  \end{proof}

\begin{proposition} \label{prop: summingbasissubseq} Suppose that $(e_i)$ is a nontrivial weakly Cauchy dictionary for $X$. If every subsequence of $(e_i)$ has the NQP for its closed linear span
(in particular, if $(e_i)$ has the CQP)
then $(e_i)$ has a subsequence equivalent to the summing basis of $c_0$. \end{proposition}
\begin{proof} By Theorem~\ref{thm: indepinterest} either $(e_i)$ has a subsequence equivalent to the summing basis or a basic subsequence $(e_{n_i})$ such that
$(e^*_{n_i})$ is weakly null in $[e_{n_i}]^*$. But $(e_{n_i})$ has the NQP for its closed linear span, so the latter alternative cannot happen by Proposition~\ref{prop: nosubseqweaklynull}. \end{proof}

\begin{proof}[Proof of Theorem~\ref{thm: conversemainthm}]
By Rosenthal's $\ell_1$ theorem \cite{R1}, either $(e_i)$ has a subsequence that is equivalent to the unit vector basis of $\ell_1$ or a weakly Cauchy basic subsequence.
The first possibility cannot occur since the unit vector basis of $\ell_1$ does not have the NQP. For the second possibility, either the  subsequence is weakly null or it is nontrivial weakly Cauchy. In the former case there is a subsequence equivalent to the unit vector basis of $c_0$ by Proposition~\ref{prop: c0subseq}, and in the latter there
is a subsequence equivalent to the summing basis by Proposition~\ref{prop: summingbasissubseq}.
\end{proof}
Combining  Theorem~\ref{thm: existence} and Theorem~\ref{thm: conversemainthm}
we obtain a new characterization of separable Banach spaces  containing $c_0$

\begin{theorem} Let $X$ be a separable Banach space. The following are equivalent: \begin{itemize}
\item[(a)] $c_0\hookrightarrow X$;
\item[(b)] $X$ has a weakly null bounded and total minimal system with the SCQP;
\item[(c)] $X$ has a total minimal system $(e_i)$ with the CQP;
\item[(d)] $X$ has a dictionary $(e_i)$ with no nonzero weak limit point such that every subsequence of $(e_i)$ has the NQP for its closed linear span.
\end{itemize} \end{theorem}

\begin{proof} (a)~$\Rightarrow$~(b) follows from Theorem~\ref{thm: existence};  (b)~$\Rightarrow$~(c) is trivial; (c)~$\Rightarrow$~(d) follows from the fact that a total minimal system has no nonzero subsequential weak limit point.
To prove (d)~$\Rightarrow$~(a), note that $(e_i)$ has a weakly Cauchy basic subsequence, so the result follows from Propositions~\ref{prop: c0subseq} and \ref{prop: summingbasissubseq}.

\end{proof}

We conclude this section with some results about NQP minimal systems that are motivated by Problem~\ref{prob: NQPcontainc0} above.

\begin{proposition} \label{prop: nosubseqnontrivialweakcauchy} \mnote{This and the next are based on faxes of Ted.}
 Let $(e_i)$ be a minimal system for $X$ with the NQP. Then no subsequence of $(e^*_i)$ is nontrivial weakly Cauchy. \end{proposition}
\begin{proof} Suppose that $(e_i)$ has the $(\varepsilon,\delta)$-NQP and that $(e^*_{n_i})$ is nontrivial weakly Cauchy. After passing to a subsequence of $(e^*_{n_i})$ we may assume that $(f_i)$
is a weakly null basis for $Y=[(f_i)] = [(e_{n_i}^*)]$ , where $f_1 = e^*_{n_1}$ and $f_i = e^*_{n_i} - e^*_{n_{i-1}}$ for $i \ge2$ \cite[Prop. II.1.7]{AGR}.
Setting $e = \sum e_{n_i}|_Y$ (the sum converging weak-star in $Y^*$), $n_0:=0$, and  $e_0:=0$,
the sequence of biorthogonal functionals $(f_i^*) \subset Y^*$ is given by $f^*_i= e - \sum_{j=0}^{i-1} e_{n_j}|_Y$. We claim that
$(f^*_i)$ has the NQP for its closed linear span in $Y^*$. To check this claim,
let $x = \sum_{i \in A} a_i f^*_i$, where $A \subseteq{N}$ is finite. Then we may rewrite the expression for $x$ in the form  \begin{equation} \label{eq: expressionforx}
x = b f_1^* + \sum_{i \in B} b_i e_{n_i}|_Y \end{equation} for some finite $B \subset \mathbb{N}$ and scalars $b, b_i$.
 Since $(e_i)$ has the $(\varepsilon,\delta)$-NQP there exists $z = \sum_{i \in C} m_i \delta e_i$ ($m_i \in \mathbb{Z}$), where $C$ is a finite subset of
$\mathbb{N}$,
such that $\|\sum_{i \in B} b_i e_{n_i} - z\| \le \varepsilon$. Since $e_{n_i}|_Y = f^*_i - f^*_{i+1}$, it follows that \begin{equation} \label{eq: expressionforz}
   z|_Y = \sum_{i \in C'} m'_i \delta f^*_i \end{equation} for some finite $C' \subset \mathbb{N}$ and  $m_i' \in \mathbb{Z}$.
Choose $m \in \mathbb{Z}$ such that $|b-m\delta| \le \delta$. From \eqref{eq: expressionforx}  and \eqref{eq: expressionforz}, we obtain
\begin{align*}
\|x - (m\delta f^*_1 + \sum_{i \in C'} m'_i \delta f^*_i)\|_{Y^*} &\le |m\delta-b|\|f^*_1\| + \|\sum_{i \in B} b_i e_{n_i}|_Y - z|_Y\|_{Y^*} \\
&\le \delta \|f^*_1\| + \|\sum_{i \in B} b_i e_{n_i} - z\| \\ &\le \delta \|f^*_1\| + \varepsilon, \end{align*}
which verifies the claim. Thus $(f^*_i)$ has the NQP for its closed linear span and its biorthogonal sequence $(f_i)$ is weakly null. But this contradicts Proposition~\ref{prop: nosubseqweaklynull}.  \end{proof}

\begin{theorem} Let $(e_i)$ be a seminormalized basis with the NQP. Then every subsequence of $(e^*_i)$ has a further subsequence equivalent to
the unit vector basis of $\ell_1$. \end{theorem}

\begin{proof} By Proposition~\ref{prop: nosubseqweaklynull} no subsequence of $(e_i^*)$ is weakly null, and by Proposition~\ref{prop: nosubseqnontrivialweakcauchy} no subsequence is nontrivial weakly Cauchy. Thus, by Rosenthal's $\ell_1$ theorem, every subsequence of $(e^*_i)$ has a further subsequence equivalent to
the unit vector basis of $\ell_1$. \end{proof}

\section{Some Notions Related to the CQP}\label{openquestions}

There seems to be very little known about the relationships
 between the different quantization properties
 introduced in the previous sections. Let us recast some of
 the questions we  formulated in previous sections.

 Throughout this section $(e_i)$ and $(e^*_i)$ is   a bounded minimal
 system of a Banach space $X$ and we assume that $(e_i)$ (and, thus, also $(e^*_i)$) are semi-normalized.
\bigskip
\begin{qst}\label{Q:oq1} Let $\vp,\delta>0$.

\begin{enumerate}
\item  If $(e_i)$ satisfies the $(\vp,\delta)$-CQP, does it satisfy the $(\vp,\delta/2)$-SCQP,
 does it satisfy $(\vp,\delta)$-neighborly CQP (see Remark~\ref{rem: neighborly})?

In the case that the answer to our aforementioned questions are negative  do they at least
 have  qualitative  positive answers, i.e. does the CQP imply the SCQP,
 does the CQP imply the neighborly CQP?

\item In our next example we will exhibit that for some $\vp,\delta>0$ the
 $(\vp,\delta)$-NQP does not imply the $(\vp,\delta/2)$-SNQP.
  But we do  not know whether or not the NQP implies the SNQP.
\end{enumerate}
\end{qst}
One can reformulate these questions into finite dimensional ones.
Assume that $n\in\N$ and that  $K\subset \R^n$ is a symmetric and convex body (i.e. $0\in {\mathop K}^\circ$).

Let us consider the following properties $K$ may have
\begin{align}
\bigcup_{z\in \Z^n} z+K&=\R^n \tag{P1}\\
\bigcup_{z\in \prod_{i=1}^n D_i} z+K&=\R^n  \tag{P2}\\
\qquad\text{ whenever } &D_i\subset R, 0\in D_i, D_i\text{ is $\frac12$-net for }i=1,2,\ldots, n
\notag\\
[0,1]^n&\subset \bigcup_{\vp=(\vp_1,\vp_2,\ldots ,\vp_n)\in\{0,1\}^n} \vp+K
\tag{P3}
\end{align}
Note that (P3) means that, not only is every point of $\R^n$ an element
of some translate of $K$ by some point $p$ having integer coordinates, but that
  $p$ can be chosen so that $\max_{i=1,2\ldots, n}|x_i-p_i|\le 1$.

It is easy to see that $(e_i)$ satisfies $(\vp,\delta)$-CQP, $(\vp,\delta/2)$-SCQP
 or $(\vp,\delta)$-neighborly CQP, if  and only if for any finite $I\subset \N$
the set
$$K_I=\frac{\vp}{\delta} B_X \cap [e_i:i\in I],$$
satisfies (P1), (P2) or (P3) respectively.

If we do not assume that $(e_i)$ is a monotone basis a similar
 statement for NQP and SNQP is slightly more complicated.

First if $E=(R^n,\|\cdot\|)$ is  finite dimensional then
 the unit vector basis $(e_i)$ has the $(\vp,\delta)$-NQP or the
 $(\vp,\delta/2)$-SNQP if and only if $\frac{\vp}{\delta}B_E$ satisfies
 (P1) or (P2).
If for all $n\in\N$
$K_{\{1,2,\ldots, n\}}$ (defined as above) satisfies (P1) or (P2) then for any
 $\eta>0$
 $(e_i)$ has the $(\vp,\delta-\eta)$-NQP or the $(\vp,\delta-\eta)$-SNQP
 respectively.
Conversely, if $(e_i)$ is a monotone basis which satisfies the $(\vp,\delta)$-NQP or the $(\vp,\delta)$-SNQP,
 then for all $n\in\N$  the set $K_{\{1,2,\ldots n\}}$  satisfies (P1) or (P2), respectively.

The following example shows that $(P1)\not\Rightarrow  (P2)$.

\noindent{\bf Example}.
In $\R^2$ let $K$ be the convex hull of the points
$$ P_1=(\frac14, 1),\,
   P_2=(\frac34, 1),\,
P_3=(-\frac14, -1),\,\text{ and }
P_4=(-\frac34,-1).
$$
Instead of a  formal proof, we leave it to the reader to verify the following by drawing a picture:
\begin{enumerate}
\item[a)] $K$ is a parallelogram which tiles $R^n$, i.e.
$$\bigcup_{z\in\Z^2} z+K=\R^2\text{ and }
(z+{\mathop K}^\circ) \cap (z'+{\mathop K}^\circ)
\text{ whenever }z\not=z'\text{ are in }\Z^2.$$
\item[b)]  For $\Q=\frac34 P_2+ \frac14 P_3$ we have
$$Q\in [(0,0)+K]\cap
\left[(1,1)+K\right].$$
\item[c)] For small enough $\eta>0$
$$P- ( 0, \eta/4)\not\in
\bigcup_{z\in \Z\times(1-\eta)\Z} z+ K.$$
(thus $K$ does not satisfy $(P2)$).
\end{enumerate}
The aformentioned questions can be now reformulated as follows.
\begin{qst}\label{Q:oq2} Let $K\subset \R^n$ be convex and symmetric and put for
 $I\subset\{1,2,\ldots, n\}$
$$K_I=\big\{(x_1,x_2,\ldots,x_n)\in K: x_i=0 \text{ for }i\in\{1,2,\ldots, n\}\setminus I\big\}.$$
\begin{enumerate}
\item If $K_I$ satisfies $(P1)$ for all $I\subset\{1,2\ldots, n\}$,
 does it satisfy $(P2)$  or $(P3)$?
\item Is there at least a universal constant $c\ge 1$ so that
 if $K_I$ satisfies $(P1)$ for all $I\subset\{1,2,\ldots, n\}$, then it satisfies
 $(P2)$  or $(P3)$?
\item Is there  a universal constant $c\ge 1$ so that
 if $K$ satisfies $(P1)$ then it satisfies $(P2)$?
\end{enumerate}

\end{qst}


\begin{thebibliography}{WW}

\bibitem{AGR}
 Spiros A. Argyros, Gilles Godefroy and Haskell P. Rosenthal, {\em Descriptive Set Theory and Banach Spaces}
in: William B. Johnson and Joram Lindenstrauss (eds.), Handbook on the Geometry of Banach Spaces Vol. 2, North Holland, Amsterdam, 2003, 1007-1069.
497-532.
\bibitem{BL} B. Beauzamy and J.-T. Laprest\'e,
{\em Mod\`eles \'etal\'es des espaces de Banach},
Travaux en cours, Hermann, Paris, 1984.

\bibitem{BP1} C. Bessaga and A. Pe\l czy\'nski, {\em On bases and unconditional
convergence of series in Banach spaces}, Studia Math. {\bf 17} (1958), 151--164.



\bibitem{BP2} C. Bessaga and A. Pe\l czy\'nski, {\em Spaces of Continuous Functions IV
(On isomorphic classification of spaces $C(S)$)}, Studia Math. {\bf 19} (1960), 53--62.

\bibitem{DOSZ} S. J. Dilworth, E. Odell, Th. Schlumprecht, and Andr\'as Zs\'ak,
{\em Partial Unconditionality}, preprint, 2005.

\bibitem{DKK} S. J. Dilworth, N. J. Kalton and Denka Kutzarova, {On the existence of
almost greedy bases in Banach spaces. Dedicated to Professor Aleksander Pe\l czy\'nski
on the occasion of his 70th birthday},
Studia Math. {\bf 159} (2003), no. 1, 67-101.


\bibitem{DKKT} S. J. Dilworth, N. J. Kalton, Denka Kutzarova and V. N. Temlyakov,
{\em The thresholding greedy algorithm, greedy bases, and duality},
Constr. Approx.  {\bf 19} (2003), 575--597.


\bibitem{E1} John Elton, {\em Weakly null normalized sequences in Banach spaces},
Ph.D. thesis, Yale University, 1978.

\bibitem{E2} John Elton, {\em Extremely weakly unconditionally convergent series},
Israel J. Math {\bf 40} (1981), 255--258.

\bibitem{GOW} I. Gasparis, E. Odell and B. Wahl, {\em Weakly null sequences in the
Banach spaces C(K)}, to appear.

\bibitem{G} W. T. Gowers, {\em Lipschitz functions on classical spaces}, European J. Combin. {\bf 13} (1992), 141--151.

\bibitem{HOR}
R. Haydon, E. Odell and H. Rosenthal, {\em On certain classes of Baire-$1$ functions with applications to Banach space theory}. Functional analysis (Austin, TX, 1987/1989), 1--35, Lecture Notes in Math., 1470, Springer, Berlin, 1991.


\bibitem{J} R. C. James, {\em Uniformly non-square Banach spaces}, Ann. of Math. {\bf 80} (1964), 542--550.


\bibitem{KT} S. V. Konyagin and V. N. Temlyakov,  {\em A remark
on greedy approximation in Banach spaces}, East J.
Approx. {\bf 5} (1999), no. 3, 365--379.







\bibitem{LT}
J.\ Lindenstrauss and L. Tzafriri, {\em Classical Banach
spaces I, Sequence Spaces}, Springer-Verlag, Berlin-Heidelberg, 1977.

\bibitem{LAT} J. Lopez-Abad and S. Todorcevic, {\em Pre-compact families of integers and weakly null sequences in Banach spaces}, preprint, 2005.

\bibitem{M} A. A. Milutin, {\em Isomorphisms of spaces of continuous functions on compacta of power continuum}, Teori Func. (Kharkov) {\bf 2} (1966), 150--156 (Russian).


\bibitem{O} E. Odell, {\em Applications of Ramsey theorems to Banach space theory}
in:
Notes in Banach spaces, (H. E. Lacey, ed.), 379--404, Univ. Texas Press, Austin, TX,
1980.
\bibitem{OS} E. Odell and Th. Schlumprecht, {\em The distortion problem}, Acta Math. {\bf 173} (1994), 259--281.


\bibitem{OP}
R.~I. Ovsepian and A. Pe\l czy{\'n}ski, {\em On the existence of a
  fundamental total and bounded biorthogonal sequence in every separable
  Banach space, and related constructions of uniformly bounded orthonormal
  systems in $L^2$}, Studia Math. \textbf{54} (1975), no.~2, 149--159.

\bibitem{P2}
A.~Pe{\l}czy{\'n}ski,
 \emph{All separable {B}anach spaces admit for every $\varepsilon >0$\
  fundamental total and bounded by $1+\varepsilon $ biorthogonal sequences},
  Studia Math. \textbf{55} (1976), no.~3, 295--304.




\bibitem{R1} H. P. Rosenthal, {\em A characterization of Banach spaces containing $\ell_1$},
Proc. Nat. Acad. Sci. {\bf 71} (1974), 2411--2413.



\bibitem{R2} H. P. Rosenthal, {\em A characterization of Banach spaces containing $c_0$},
J. Amer. Math. Soc. {\bf 7} (1994), 707--748.




\bibitem{S} A. Sobczyk, {\em Projection of the space $(m)$ on its subspace $(c_0)$},
Bull. Amer. Math. Soc. {\bf 47} (1941), 938--947.


\bibitem{W2}
P. Wojtaszczyk, {\em Existence of some special bases in Banach
spaces}, Studia Math. {\bf 67} (1973), 83--93.

\bibitem{W1} P. Wojtaszczyk, {\em Greedy algorithm for general biorthogonal
systems}, J. Approx. Theory {\bf 107} (2000), 293--314.

\end{thebibliography}
\end{document}